\documentclass[11pt]{amsart}
\usepackage[margin=1.15in]{geometry}
\usepackage{amsmath,amssymb,amsthm,mathtools}
\usepackage{enumitem}
\usepackage{booktabs,array}
\usepackage{microtype}
\usepackage{fancyvrb}
\usepackage{tikz}
\usepackage{amsaddr}
\usetikzlibrary{arrows.meta,positioning}
\usepackage[colorlinks=true,linkcolor=blue!50!black,citecolor=green!35!black,urlcolor=blue!50!black]{hyperref}

\newtheorem{mainthm}{Theorem}
\newtheorem{mainprop}[mainthm]{Proposition}
\newtheorem{theorem}{Theorem}[section]
\newtheorem{proposition}[theorem]{Proposition}
\newtheorem{lemma}[theorem]{Lemma}
\newtheorem{corollary}[theorem]{Corollary}

\theoremstyle{definition}

\theoremstyle{remark}
\newtheorem{remark}[theorem]{Remark}
\numberwithin{equation}{section}

\newcommand{\E}{\mathbb{E}}
\newcommand{\EW}{\mathbb{E}_W}
\newcommand{\1}{\mathbf{1}}
\DeclareMathOperator{\re}{Re}
\DeclareMathOperator{\im}{Im}
\DeclareMathOperator{\Leg}{Leg}
\DeclareMathOperator{\tr}{tr}
\newcommand{\gap}{\mathrm{gap}}
\newcommand{\Ksq}{\mathcal K}

\title[Simple zeros of the zeta-function]{Simple zeros of the Riemann zeta-function:\\ a local multiplicity test and long mixed moments}
\author{Yuanda Xu}
\email{yuanda@math.princeton.edu}
\date{}

\begin{document}

\begin{abstract}
Assuming the Riemann Hypothesis, we show that at least $71.274066\%$ of the nontrivial zeros of $\zeta(s)$ are simple, and at least $85.368230\%$ are distinct. The previous bounds were $19/27\approx70.37\%$, due to Conrey, Ghosh and Gonek and to Bui and Heath-Brown (2013), and $84.77\%$, due to Chirre, Gon\c{c}alves and de Laat (2020). The method combines the two classical approaches to simple zeros. Montgomery's pair correlation argument rests on a pointwise consequence of the explicit formula: at a zero of multiplicity $m\ge2$, a Dirichlet polynomial $J$ over prime powers of length $T^{\psi}$, $\psi<1$, is forced to have real part at most $-a<0$. The mollified derivative $\zeta'(\rho)M(\rho)$ of Conrey, Ghosh and Gonek vanishes there. Both facts are used at the same zero, in a single least-squares inequality. The new analytic input is the mixed moment $\sum_\rho\zeta'(\rho)M(\rho)J(1-\rho)$ with $J$ much longer than the mollifier. We evaluate it under the Riemann Hypothesis alone; the long prime sum stays in the numerator and the mollifier variable is the modulus, so the barrier at $T^{1/2}$ does not arise. A second-order version of the argument gives $0.70736$. We show that no second-order certificate of this type can exceed $0.7073646063$, and pass this ceiling with fourth moments of short Dirichlet polynomials and with mixed moments against products of a bounded number of primes. All constants are certified in exact rational arithmetic.
\end{abstract}

\maketitle

\section{Introduction}\label{sec:intro}

Let $N(T)$ be the number of zeros $\rho=\beta+i\gamma$ of $\zeta(s)$ with $0<\gamma\le T$, counted with multiplicity, let $N_s(T)$ be the number of simple zeros in the same range, and let $N_d(T)$ be the number of distinct zeros. It is conjectured that every zero is simple. Conrey~\cite{Con89} proved unconditionally that more than two fifths of the zeros are simple and on the critical line, and recent work has raised this to more than two thirds~\cite{AF26,Lam26,Kna26}. The subject of this paper is the bound for $N_s(T)/N(T)$ under the Riemann Hypothesis (RH).

\subsection{Two routes to simple zeros}
There are two classical methods, and until now they have been developed separately.

\emph{Pair correlation.} Montgomery~\cite{Mon73} studied the zeros through the function
\[
  F(\alpha)=\frac{2\pi}{T\log T}\sum_{0<\gamma,\gamma'\le T}T^{i\alpha(\gamma-\gamma')}w(\gamma-\gamma'),\qquad w(u)=\frac4{4+u^2}.
\]
The explicit formula expresses $F(\alpha)$ for $|\alpha|<1$ in terms of primes, and Montgomery computed it there, assuming RH. Since $F$ is a Fourier transform of the pair statistic, one can test it against a kernel $K\ge0$ with $\widehat K$ supported in $[-1,1]$. The diagonal terms $\gamma=\gamma'$ include each zero $m$ times, so $\sum_\rho m_\rho$ is bounded above by an average of $K$ over pairs of zeros. This gave $N_s(T)\ge(\frac23+o(1))N(T)$. Optimising the kernel improves this to $0.6725$ (Montgomery and Taylor; see~\cite{Mon74,CG93}), and to $0.6727$ (Cheer and Goldston~\cite{CG93}). Chirre, Gon\c{c}alves and de Laat~\cite[Corollary~2]{CGdL20} used semidefinite programming to obtain $0.6792$ under RH, and $0.6845$ assuming in addition the generalised Riemann Hypothesis (GRH) for Dirichlet $L$-functions.

\emph{Mollified discrete moments.} Conrey, Ghosh and Gonek~\cite{CGG98} looked at the values of $\zeta'$ at the zeros. Let $M(s)=\sum_{h\le T^\theta}\mu(h)P(\cdots)h^{-s}$ be a mollifier. The product $\zeta'(\rho)M(\rho)$ vanishes at every multiple zero, so by Cauchy--Schwarz
\[
  \Bigl|\sum_\rho\zeta'(\rho)M(\rho)\Bigr|^2\le N_s(T)\sum_\rho|\zeta'(\rho)M(\rho)|^2 .
\]
Both sums can be evaluated for $\theta<\frac12$. Conrey, Ghosh and Gonek obtained $N_s(T)\ge(\frac{19}{27}+o(1))N(T)$, assuming RH and the generalised Lindel\"of hypothesis. Bui and Heath-Brown~\cite{BHB13} removed the Lindel\"of hypothesis; they proved $19/27=0.7037\ldots$ under RH alone. This bound has stood since 2013. Longer mollifiers would give more, but the mean values are not known for $\theta\ge\frac12$.

The two methods use quite different information. The pair correlation route uses primes up to $T^{1-\varepsilon}$, but only through averages over all zeros. The mollifier route works zero by zero, but its Dirichlet polynomials are shorter than $T^{1/2}$.

\subsection{Results}
Alp\"oge has announced a proof of the Riemann Hypothesis~\cite{Alp26}.
We state all results assuming RH; nothing else is assumed in Theorems~\ref{thm:main} and~\ref{thm:distinct}.

\begin{mainthm}\label{thm:main}
Assume RH. Then
\[
  \liminf_{T\to\infty}\frac{N_s(T)}{N(T)}\ \ge\ 0.7127406684 .
\]
\end{mainthm}

\begin{mainthm}\label{thm:distinct}
Assume RH. Then
\[
  \liminf_{T\to\infty}\frac{N_d(T)}{N(T)}\ \ge\ 0.8536823046 .
\]
\end{mainthm}

The constant in Theorem~\ref{thm:main} is a truncation of the exact value $0.712740668418879614\ldots$ of a rational dual certificate (Section~\ref{sec:higher} and Appendix~\ref{app:data}). Theorem~\ref{thm:distinct} does not follow from Theorem~\ref{thm:main}: the implication $N_d/N\ge(1+N_s/N)/2$ is false in general. It comes from a separate certificate that controls zeros of multiplicity three or more and the packing of double zeros (Section~\ref{sec:distinct}). The previous best bound under RH was $0.8477$, obtained by Chirre, Gon\c{c}alves and de Laat~\cite[Corollary~3]{CGdL20} from $19/27$ and their pair correlation bound for $\sum_\rho m_\rho$; assuming GRH as well, they obtain $0.8486$.

The proof has three steps, and the intermediate bounds show where the gain comes from.
\begin{center}
\begin{tabular}{@{}l>{\raggedright\arraybackslash}p{5.6cm}l@{}}
\toprule
Certificate & Moments used & Bound\\
\midrule
Bui--Heath-Brown~\cite{BHB13} & second moments of $\zeta'M$ & $19/27=0.7037037\ldots$\\
Theorem~\ref{thm:second} & $+$ second moments of $J$, mixed moment $\zeta'M\cdot J$ & $0.7073645768\ldots$\\
Theorem~\ref{thm:fourth} & $+$ fourth moments of short polynomials & $0.7094399755\ldots$\\
Theorem~\ref{thm:main} & $+$ mixed moments against long products of primes & $0.7127406684\ldots$\\
\bottomrule
\end{tabular}
\end{center}
The second-order certificate cannot be pushed much further. Over all kernels and all mollifier tapers it is bounded by $0.7073646063$ (Proposition~\ref{prop:ceiling}). So higher moments are needed to go beyond it, and they give most of the improvement.

The packing argument behind Theorem~\ref{thm:distinct} also applies to families. For primitive Dirichlet $L$-functions, with the weights of Chandee, Lee, Liu and Radziwi\l\l~\cite{CLLR14}, Chirre, Gon\c{c}alves and de Laat~\cite[Theorem~5]{CGdL20} obtained $0.9350$ for simple zeros under GRH. Adding packing to the pair correlation bound gives a small improvement.

\begin{mainprop}\label{prop:family}
Assume GRH. In the setting of~\cite{CLLR14} (made precise in Section~\ref{sec:family}), the weighted proportion of simple zeros of primitive Dirichlet $L$-functions is at least $0.93508798998$, and the weighted proportion of distinct zeros is at least $0.96754399499$.
\end{mainprop}

\subsection{The idea}\label{sec:idea}
Montgomery's argument rests on a pointwise inequality, which is then averaged over all zeros. We use it without averaging.

\begin{figure}[ht]
\centering
\begin{tikzpicture}[x=1.25cm,y=2.1cm]
  \draw[->] (-3.4,0) -- (3.6,0) node[right] {$x=(\gamma'-\gamma)L/2\pi$};
  \draw[->] (0,-0.05) -- (0,1.25);
  \draw[thick,blue!60!black,domain=-3.3:3.3,samples=199,smooth] plot ({\x},{(sin(180*\x)/(3.14159265*\x))^2});
  \node[blue!60!black] at (1.9,0.75) {$K(x)$};
  \foreach \p in {-2.55,-1.35,0.72,1.61,2.9} {\fill[red!70!black] (\p,0) circle (1.6pt);}
  \fill[red!70!black] (0,0) circle (2.4pt);
  \node[below] at (0,-0.04) {\small $\gamma$ (double)};
  \draw[dashed] (-3.3,1) -- (3.3,1) node[right] {\small $1$};
  \node[align=left,anchor=west] at (3.75,0.45) {\small sum over all zeros $\ge m_\rho$};
\end{tikzpicture}
\caption{The Fej\'er kernel $K(x)=(\sin\pi x/\pi x)^2$ centred at a zero $\gamma$, with neighbouring zeros. If $\gamma$ is a double zero, the sum of $K$ over all zeros, with multiplicity, is at least $2$. The explicit formula writes this sum as $k(0)-2\re J(\rho)+o(1)$.}\label{fig:kernel}
\end{figure}
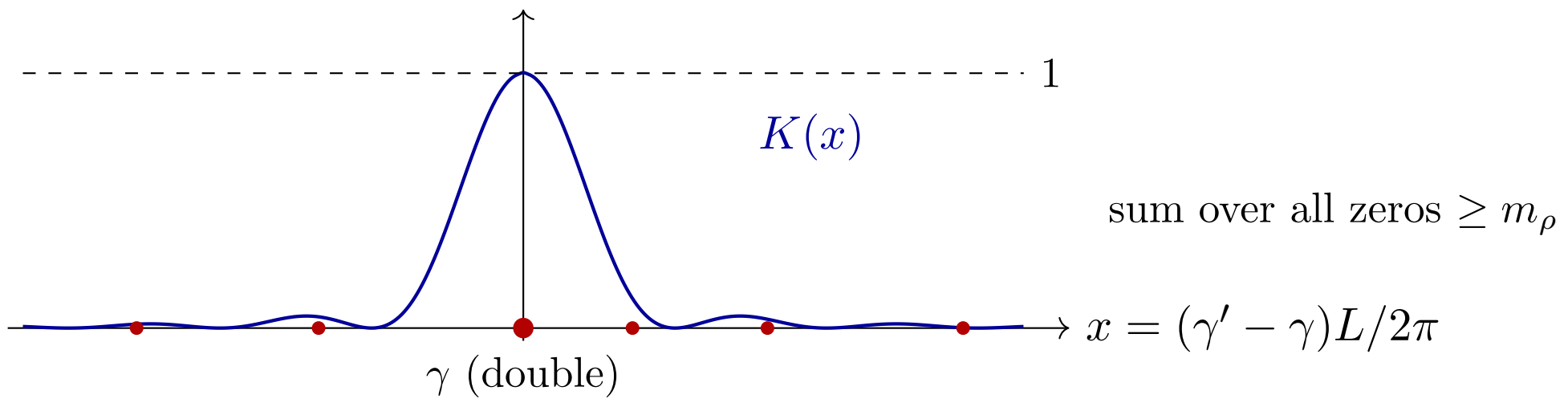

Let $K\ge0$ be an even function on $\mathbb R$ with $K(0)=1$, whose Fourier transform $k=\widehat K$ is supported in $[-\psi,\psi]$ with $\psi<1$. Put $L=\log(T/2\pi)$. Fix a zero $\rho=\frac12+i\gamma$ of multiplicity $m_\rho$, with $T\le\gamma\le2T$. Every term in $\sum_{\gamma'}K((\gamma'-\gamma)L/2\pi)$ is nonnegative, and the zero $\rho$ itself contributes $m_\rho$. The explicit formula evaluates the sum (Lemma~\ref{lem:local}), and we get
\begin{equation}\label{eq:intro-local}
  m_\rho\ \le\ \sum_{\gamma'}K\Bigl(\frac{(\gamma'-\gamma)L}{2\pi}\Bigr)=k(0)-2\re J(\rho)+O\Bigl(\frac1L\Bigr),\qquad J(s)=\frac1L\sum_{n}\frac{\Lambda(n)k(\log n/L)}{n^{s}} .
\end{equation}
The polynomial $J$ has length $(T/2\pi)^{\psi}$. If $\rho$ is a multiple zero, the left side is at least $2$, and so
\begin{equation}\label{eq:intro-sign}
  \re J(\rho)\le-a+O(L^{-1}),\qquad a:=1-\frac{k(0)}2 .
\end{equation}
For $\psi$ close to $1$ one can arrange $a>0$; for the Fej\'er kernel with $\psi=1$, $a=\frac12$. So the prime polynomial $J$ \emph{sees} multiplicity, zero by zero: at every multiple zero it is pushed to the left of $-a$. Figure~\ref{fig:kernel} illustrates~\eqref{eq:intro-local}.

Now let $F(\rho)=\zeta'(\rho)M(\rho)/L$ be a mollified detector. It vanishes at every multiple zero. Put $q=J+a$ and fix $t\ge0$. At a multiple zero, $F=0$ and $\re q\le o(1)$, so
\[
  \re\bigl(1-F(\rho)-t\,q(\rho)\bigr)\ \ge\ 1-o(1).
\]
Hence the indicator of the multiple zeros is bounded by $|1-F-tq|^2$, up to $o(1)$, and
\begin{equation}\label{eq:intro-cert}
  1-\frac{N_s(T)}{N(T)}\ \le\ \frac1{N(T)}\sum_{0<\gamma\le T}\bigl|1-F(\rho)-t\,q(\rho)\bigr|^2+o(1).
\end{equation}
The right side is a quadratic polynomial in $t$ and in the coefficients of the mollifier. Its coefficients are the first and second moments of $F$ and $J$ over the zeros. The moments of $F$ alone are those of Conrey, Ghosh and Gonek and Bui and Heath-Brown. The moments of $J$ alone are classical; they follow from Landau's formula, in Gonek's uniform version~\cite{Gon93}. The only new quantity is the mixed moment
\[
  \frac1{N(T)}\sum_{0<\gamma\le T}F(\rho)\overline{J(\rho)} .
\]
With $t=0$, minimising~\eqref{eq:intro-cert} over the mollifier gives back $19/27$. With $F=0$, it gives $(\mu+a)^2/(\E|J+a|^2)$, a Montgomery-type bound, which for the Fej\'er kernel is exactly Montgomery's $2/3$. Together they give more than either. The reason is that the two tests are correlated but not identical. The mollifier predicts part of $J$, through their covariance, and the part of $J+a$ that it does not predict still has positive mean. This residual mean is the gain.

\subsection{A proof sketch: the simplest certificate}\label{sec:sketch}
Here is the whole argument with the smallest data: the Fej\'er kernel, one mollifier, and $(\theta,\psi)$ at the endpoint $(\frac12,1)$. All limits are as $T\to\infty$, and the endpoint is reached as a limit of admissible exponents (Section~\ref{sec:second}). Write $\E$ for the average over the zeros $0<\gamma\le T$.

\emph{The kernel.} $K(x)=(\sin\pi x/\pi x)^2$, $k(u)=(1-|u|)_+$, $k(0)=1$, $a=\frac12$.

\emph{The detector.} $F=\zeta'(\rho)M(\rho)/L$ with the Bui--Heath-Brown mollifier, whose taper is $f(u)=1-u-2u^2$ in the variable $u=\log h/L\in[0,\frac12]$. Bui and Heath-Brown's theorem gives
\[
  \E F\to A=\tfrac{19}{24},\qquad \E|F|^2-|\E F|^2\to V=\tfrac{19}{72} .
\]
\emph{The witness.} By Landau's formula, $\E J\to\mu=-\int_0^1u\,k(u)\,du=-\frac16$ and $\E|J|^2\to\int_0^1u\,k(u)^2du=\frac1{12}$. So the variance is $V_J=\frac1{18}$, and $\E(J+a)\to\frac13$.

\emph{The mixed moment.} Theorem~\ref{thm:longmixed} below gives the real part of the covariance of $F$ and $J$:
\[
  C=\int_0^1k(u)\,\mathfrak c(u)\,du=\tfrac{203}{1920},\qquad
  \mathfrak c(u)=\begin{cases}u-\frac23u^3,&0\le u\le\frac12,\\ \frac{19}{24}-\frac u2-\frac{u^2}2,&\frac12\le u\le1.\end{cases}
\]
\emph{The least-squares problem.} Minimise $\E|1-\beta F-t(J+a)|^2$ over real $\beta$ and $t\ge0$. With the means $m=(A,\frac13)$ and the covariance $\Sigma=\begin{psmallmatrix}V&C\\C&V_J\end{psmallmatrix}$, the minimum is $1-\Ksq/(1+\Ksq)$ with $\Ksq=m^{\mathsf T}\Sigma^{-1}m$ (Proposition~\ref{prop:optim}). Here
\[
  \Ksq=\frac{A^2}V+\frac{\gap^2}{D_J},\qquad \gap=\frac13-\frac{AC}{V}=\frac{31}{1920},\qquad D_J=V_J-\frac{C^2}V=\frac{115519}{8755200}.
\]
The first term alone gives $\Ksq/(1+\Ksq)=19/27$. The second term is the contribution of the witness. Since $\gap>0$, the optimal $t=141360/392159$ is positive, and
\[
  \liminf_{T\to\infty}\frac{N_s(T)}{N(T)}\ \ge\ \frac{276640}{392159}=0.70542815\ldots>\frac{19}{27}.
\]
The only inputs are the inequality~\eqref{eq:intro-sign}, the moments of Bui and Heath-Brown, Landau's formula, and the mixed moment $C$. Replacing the Fej\'er kernel and the single mollifier by an optimised kernel and twelve tapers gives Theorem~\ref{thm:second}, $0.7073645768$.

\subsection{The mixed moment and the \texorpdfstring{$T^{1/2}$}{T 1/2} barrier}\label{sec:intro-mixed}
The mixed moment pairs a mollifier of length $T^\theta$, $\theta<\frac12$, with a prime polynomial of length $T^\psi$, $\psi$ close to $1$. Such long polynomials are usually out of reach. In the moments of Conrey, Ghosh and Gonek the mollifier must be shorter than $T^{1/2}$. After Gonek's lemma one meets exponential sums $e(nm/h)$ over the mollifier variables, and the products of the mollifier coefficients occur as moduli; their total length is what limits $\theta$.

The mixed moment is different. In $\sum_\rho\zeta'(\rho)M(\rho)J(1-\rho)$, Gonek's lemma produces the phase $e(pm/h)$, where $p$ is the prime of $J$, $h\le T^\theta$ is the mollifier variable and $m$ is the dual variable of $\zeta'$. The modulus is the short variable $h$, and the long prime stays in the numerator. Expanding in characters modulo $h$, we need cancellation in $\sum_p\chi(p)p^{-s}$ for moduli $h\le T^\theta$. Under GRH this is immediate. Under RH alone we use the multiplicative large sieve for large conductors and the Siegel--Walfisz theorem for small ones. This is the same substitution that Bui and Heath-Brown made for the moments of $\zeta'M$, but it is easier here, because the numerator contains a long prime. On the right side of the contour there is no diagonal at all when $J$ is a prime polynomial. For products of several primes there is one, and it has to be kept (Theorem~\ref{thm:longmixed}).

The barrier at $\theta=\frac12$ is not crossed: the mollifier is still shorter than $T^{1/2}$. What is new is that one of the two polynomials in a mixed moment may be much longer than $T^{1/2}$, as long as its long variables sit in the numerator.

\subsection{Beyond second moments}
Inequality~\eqref{eq:intro-cert} uses first and second moments only. In Section~\ref{sec:remarks} we show that, given these moments, the bound of Proposition~\ref{prop:optim} is optimal; the data are realised by a probability measure. Over all kernels and tapers, the second-order bound is at most $0.7073646063$ (Proposition~\ref{prop:ceiling}).

To go further we use two more kinds of information. The first is the support certificate of our earlier work on this problem: if $U$ is a short Dirichlet polynomial, then $\E\1_{\mathrm s}|U|^2\le\E\1_{\mathrm s}+\E(|U|^2-1)_+$, and the positive part can be bounded by fourth moments of short polynomials, which are computable under RH. We combine it with the witness. The long polynomial $J$ must not be squared, because its fourth moments are not known. So we cancel it first, using its sign at multiple zeros, and take the positive part only of a quadratic form in the short polynomials (Lemma~\ref{lem:combined}). The resulting semidefinite program gives Theorem~\ref{thm:fourth}, $0.7094399755$.

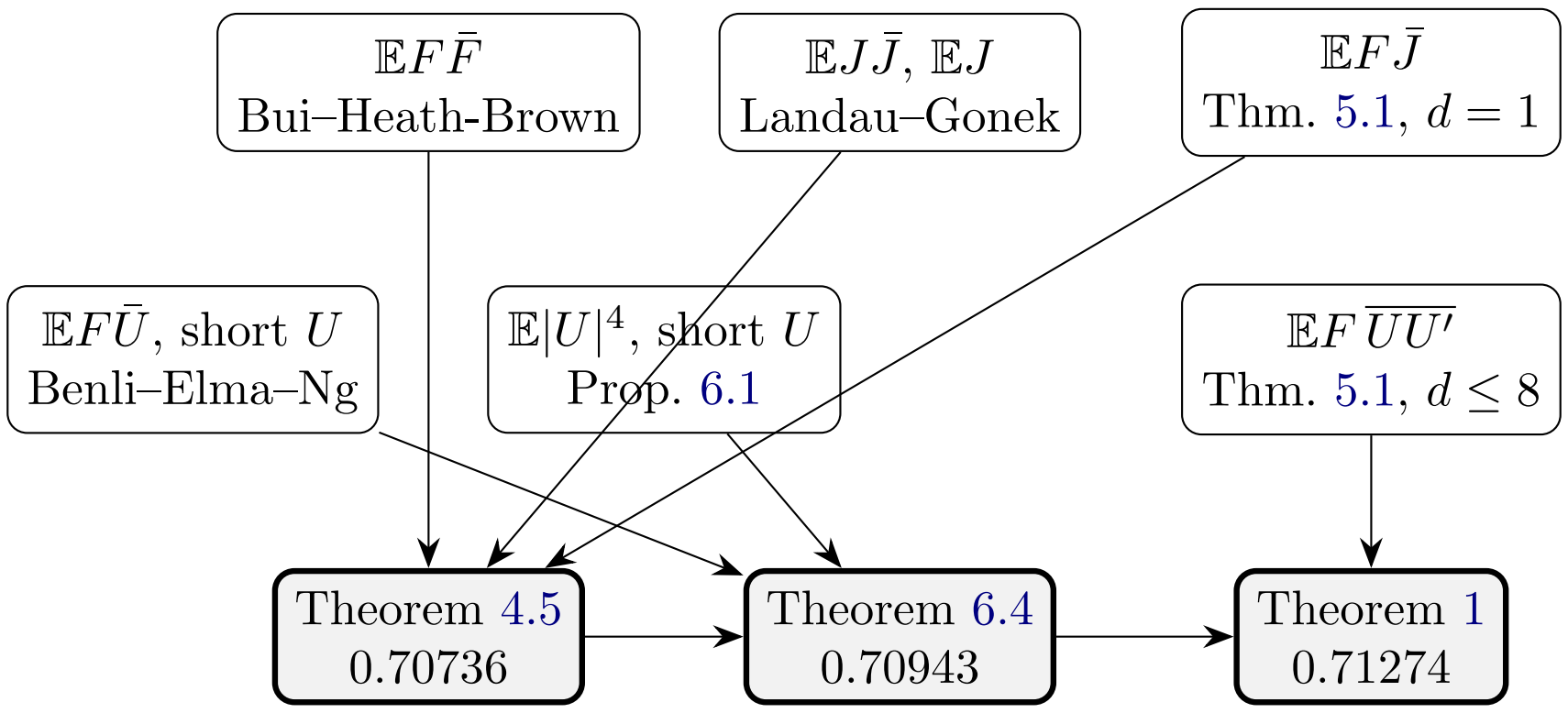
\begin{figure}[ht]
\centering
\begin{tikzpicture}[
  box/.style={draw,rounded corners,align=center,font=\small,inner sep=4pt,minimum height=0.9cm},
  res/.style={draw,very thick,rounded corners,align=center,font=\small,fill=gray!10,inner sep=4pt},
  >={Stealth[length=2.2mm]}]
  \node[box] (bhb) at (0,0) {$\E F\bar F$\\ Bui--Heath-Brown};
  \node[box] (lg) at (3.4,0) {$\E J\bar J$, $\E J$\\ Landau--Gonek};
  \node[box] (mix) at (6.8,0) {$\E F\bar J$\\ Thm.~\ref{thm:longmixed}, $d=1$};
  \node[box] (four) at (1.7,-2.0) {$\E|U|^4$, short $U$\\ Prop.~\ref{prop:pairing}};
  \node[box] (ben) at (-1.7,-2.0) {$\E F\bar U$, short $U$\\ Benli--Elma--Ng};
  \node[box] (prod) at (6.8,-2.0) {$\E F\,\overline{UU'}$\\ Thm.~\ref{thm:longmixed}, $d\le8$};
  \node[res] (r1) at (0,-4.0) {Theorem~\ref{thm:second}\\ $0.70736$};
  \node[res] (r2) at (3.4,-4.0) {Theorem~\ref{thm:fourth}\\ $0.70943$};
  \node[res] (r3) at (6.8,-4.0) {Theorem~\ref{thm:main}\\ $0.71274$};
  \draw[->] (bhb) -- (r1); \draw[->] (lg) -- (r1); \draw[->] (mix) -- (r1);
  \draw[->] (r1) -- (r2); \draw[->] (four) -- (r2); \draw[->] (ben) -- (r2);
  \draw[->] (r2) -- (r3); \draw[->] (prod) -- (r3);
\end{tikzpicture}
\caption{Which moments enter which certificate. Each certificate uses all the moments of the previous one. Here $U$ denotes a short Dirichlet polynomial, of length $T^\theta$ with $\theta<\frac12$, over products of a bounded number of primes.}\label{fig:moments}
\end{figure}

The second is the mixed moment of the detector against products of two short polynomials. Such products are supported on integers with up to eight prime factors and have length up to $T^{2\theta}$, beyond $T^{1/2}$. Theorem~\ref{thm:longmixed} evaluates these moments under RH. With them we obtain Theorem~\ref{thm:main}. Figure~\ref{fig:moments} shows which moments enter which certificate.

\subsection{What is new, and what we use}
We use the following results from the literature. From Montgomery~\cite{Mon73}, and in the form of Goldston and Gonek~\cite{GG07}, we take the explicit formula evaluated at a single zero. From Conrey, Ghosh and Gonek~\cite{CGG98} and Bui and Heath-Brown~\cite{BHB13} we take the mean values of $\zeta'M$ and $|\zeta'M|^2$ over the zeros, under RH. From Gonek~\cite{Gon84,Gon93} we take Landau's formula with a uniform error term, and the stationary phase lemma behind all discrete moments. From Benli, Elma and Ng~\cite{BEN24} we take their unconditional formula for $\sum_\rho\zeta'(\rho)X(\rho)Y(1-\rho)$ with two polynomials shorter than $T^{1/2}$.

Three things are new. First, the pointwise use of the explicit formula~\eqref{eq:intro-sign}, at the same zero as the mollified detector, inside one quadratic certificate. Second, the evaluation of the mixed moment of the detector against long polynomials over products of a bounded number of primes, under RH alone, and the observation that it avoids the barrier at $T^{1/2}$. Third, the higher-moment certificates: the cancellation of the long polynomial before the positive part is taken, the semidefinite relaxation, and the packing argument for distinct zeros.

All numerical constants come from finite-dimensional optimisation problems with rational data, and Appendix~\ref{app:data} lists the parameters needed to reproduce them.

\subsection{Organisation}
Section~\ref{sec:notation} collects notation and the results we use. Section~\ref{sec:local} proves the local multiplicity inequality~\eqref{eq:intro-local}. Section~\ref{sec:second} proves the second-order certificate, Theorem~\ref{thm:second}, assuming the mixed moment. Section~\ref{sec:longmixed} evaluates the mixed moment against long numerators, under RH. Section~\ref{sec:higher} adds fourth moments and proves Theorem~\ref{thm:main}. Section~\ref{sec:distinct} proves Theorem~\ref{thm:distinct}. Section~\ref{sec:remarks} discusses the limits of the method, open problems and the Dirichlet family. Appendix~\ref{app:data} contains the certificate data.

\section{Notation and the results we use}\label{sec:notation}

\subsection{Averages over zeros}
Throughout, $T\to\infty$ and $L:=\log(T/2\pi)$. Sums over zeros count multiplicity. We assume RH from now on, except in Section~\ref{sec:family}, where we assume GRH. Under RH, $1-\rho=\bar\rho$, so for a Dirichlet polynomial $Y(s)=\sum y(n)n^{-s}$ with real coefficients
\begin{equation}\label{eq:conj}
  Y(1-\rho)=\overline{Y(\rho)} .
\end{equation}

The mixed moment of Section~\ref{sec:longmixed} is easier to handle with a smooth weight in $\gamma$. Fix a real $w\in C_c^\infty(1,2)$, $w\not\equiv0$, put $W=w^2$, and define
\[
  N_W:=\sum_{\gamma>0}W\Bigl(\frac\gamma T\Bigr),\qquad \EW G:=\frac1{N_W}\sum_{\gamma>0}W\Bigl(\frac\gamma T\Bigr)G(\rho)
\]
for any function $G$ on the zeros. By the Riemann--von Mangoldt formula, $N_W=\frac{T}{2\pi}L\int W\,(1+O(L^{-1}))$. $\EW$ is the expectation of a probability measure on the zeros, so Cauchy--Schwarz and the usual identities for variances apply to it. Let $\1_{\mathrm s}$ be the indicator function of the simple zeros, and $m_\rho$ the multiplicity of $\rho$.

\begin{lemma}\label{lem:unweight}
Let $\phi$ be a function on the zeros with $0\le\phi\le1$ that depends only on $m_\rho$. If $\liminf_{T\to\infty}\EW\phi\ge\kappa$ for every $W$ as above, then $\liminf_{T\to\infty}N(T)^{-1}\sum_{0<\gamma\le T}\phi(\rho)\ge\kappa$.
\end{lemma}

\begin{proof}
Let $\delta>0$ and choose $W\le\1_{[1,2]}$ with $\int(\1_{[1,2]}-W)<\delta$. Since $0\le\phi\le1$ and $N(2T)-N(T)-N_W\le\frac{T}{2\pi}L\delta+O(T)$, the sum of $\phi$ over $T<\gamma\le2T$ is at least $(\kappa-2\delta-o(1))(N(2T)-N(T))$. Sum this over the intervals $(T/2^{j},T/2^{j-1}]$, $1\le j\le M$, and use $N(T/2^M)\le2^{-M}N(T)(1+o(1))$. Let $T\to\infty$, then $M\to\infty$ and $\delta\to0$.
\end{proof}

We apply the lemma with $\phi=\1_{\mathrm s}$ and with $\phi(\rho)=1/m_\rho$; the sum of $1/m_\rho$ over the zeros is $N_d(T)$. So it suffices to prove all our bounds for $\EW$.

\subsection{The mollifier}
Fix $0<\theta<\frac12$. For a real polynomial $P$ with $P(0)=0$ put
\[
  f_P(u)=\begin{cases}P(1-u/\theta),&0\le u\le\theta,\\0,&u>\theta,\end{cases}\qquad
  B_P(s)=\sum_{h\le e^{\theta L}}\frac{\mu(h)f_P(\log h/L)}{h^{s}},
\]
and
\[
  F_P(\rho)=\frac{\zeta'(\rho)B_P(\rho)}{L}.
\]
We call $f_P$ the \emph{taper} and $F_P$ a \emph{detector}. Since $P(0)=0$, $f_P$ is continuous. Every detector vanishes at every multiple zero. Write $I_f(x)=\int_x^\theta f$ and $I_f=I_f(0)$. When $P(1)=1$, $B_P$ is the mollifier of~\cite[(2)]{BHB13} with $y=e^{\theta L}$.

\begin{theorem}[Conrey--Ghosh--Gonek; Bui--Heath-Brown]\label{thm:BHB}
Assume RH, and let $0<\theta<\frac12$ and $P(0)=0$, $P(1)=1$. Then
\begin{align*}
  \sum_{0<\gamma\le T}\zeta'(\rho)B_P(\rho)&\sim\frac{TL^2}{2\pi}\Bigl(\frac12+\theta\int_0^1P\Bigr),\\
  \sum_{0<\gamma\le T}|\zeta'(\rho)B_P(\rho)|^2&\sim\frac{TL^3}{2\pi}\Bigl(\frac13+\theta\int_0^1P+\theta^2\Bigl(\int_0^1P\Bigr)^2+\frac1{12\theta}\int_0^1P'^2\Bigr).
\end{align*}
\end{theorem}

These are the two displays before~\cite[(10)]{BHB13}; the main terms go back to~\cite{CGG98}, where the Lindel\"of hypothesis for Dirichlet $L$-functions was also assumed. For $P(x)=x+\theta x(1-x)$ and $\theta\to\frac12$, the ratio of the square of the first sum to $N(T)$ times the second tends to $19/27$.

We need these moments for several tapers at once, for tapers with $P(1)\ne1$, and with the weight $W$.

\begin{lemma}\label{lem:BHBW}
Assume RH. For real polynomials $P,R$ vanishing at $0$, with $f=f_P$,
\[
  \EW F_P\to A_P:=\frac{f(0)}2+I_f,\qquad
  \re\EW F_P\overline{F_R}\to A_PA_R+V_{PR},
\]
where
\[
  V_{PR}:=\frac1{12}\Bigl(P(1)R(1)+\frac1\theta\int_0^1P'R'\Bigr).
\]
\end{lemma}

\begin{proof}
For $P=R$, $P(1)=1$ and the unweighted average, this is Theorem~\ref{thm:BHB}, since $\frac13+\theta\int P+\theta^2(\int P)^2=A_P^2+\frac1{12}$ and $\overline{F_P(\rho)}=\zeta'(1-\rho)B_P(1-\rho)/L$ by~\eqref{eq:conj}. Both sides are homogeneous of degree two in $P$, so the case $P(1)\ne0$ follows. For general $P$, the mean of $|F_{x+\varepsilon P}|^2$ is a quadratic polynomial in $\varepsilon$ and converges for every small $\varepsilon\ne0$; hence each of its coefficients converges, and polarisation gives the real part of the cross moment. The weight $W$ is introduced by partial summation over $\gamma\in[T,2T]$. For this we use that the asymptotics of~\cite{BHB13} hold uniformly for $\theta$ in compact subsets of $(0,\frac12)$ and for polynomials of bounded degree and coefficients, and that their main terms depend continuously on $\theta$; at height $t\in[T,2T]$ the same mollifier has exponent $\theta+O(L^{-1})$.
\end{proof}

At the endpoint $\theta=\frac12$ a convenient basis of tapers is
\begin{equation}\label{eq:basis}
  b_j(u)=\int_{1/2}^{u}\Leg_j(4z-1)\,dz\qquad(0\le u\le\tfrac12,\ j\ge0),
\end{equation}
where $\Leg_j$ is the Legendre polynomial. Then $b_j=f_{P_j}$ at $\theta=\frac12$ with $P_j(x)=b_j((1-x)/2)$, and $b_j(\frac12)=0$. Orthogonality of the $\Leg_j$ gives, at $\theta=\frac12$,
\begin{equation}\label{eq:basisAV}
  A_{P_0}=-\tfrac38,\quad A_{P_1}=-\tfrac1{24},\quad A_{P_j}=0\ (j\ge2),\qquad V_{P_iP_j}=\frac{\delta_{ij}}{24(2j+1)}+\frac{\delta_{i0}\delta_{j0}}{48}.
\end{equation}
For $\theta<\frac12$ we use the same $P_j$; their moments are rational functions of $\theta$ and tend to~\eqref{eq:basisAV} as $\theta\to\frac12$.

\subsection{Landau's formula and mean values of Dirichlet polynomials}
For $X>1$ let $\langle X\rangle$ be the distance from $X$ to the nearest prime power other than $X$, and put $\Lambda(X)=0$ if $X$ is not a prime power.

\begin{theorem}[Gonek~\cite{Gon85,Gon93}; see {\cite[Theorem~2]{DHPC26}}]\label{thm:LG}
Uniformly for $X,T>1$,
\begin{multline*}
  \sum_{0<\gamma\le T}X^{\rho}=-\frac{T}{2\pi}\Lambda(X)+O\bigl(X\log(2XT)\log\log(3X)\bigr)\\+O\Bigl(\log X\min\Bigl(T,\frac{X}{\langle X\rangle}\Bigr)\Bigr)+O\Bigl(\log(2T)\min\Bigl(T,\frac1{\log X}\Bigr)\Bigr).
\end{multline*}
\end{theorem}

Applied term by term to $|\mathcal A(\rho)|^2$, Theorem~\ref{thm:LG} gives an error $O(M^2\log^2T)$ for a polynomial of length $M$, which is acceptable only for $M<T^{1/2}$. We need lengths up to $T^{1-\varepsilon}$, and the following lemma shows that the error terms, summed over all pairs, are much smaller.

\begin{lemma}\label{lem:errorsum}
Assume RH. Let $2\le M\le T$, and let $a_n,c_n$ be real numbers supported on $n\le M$ with $|a_n|,|c_n|\le K$. Put $\mathcal A(s)=\sum a_nn^{-s}$ and $\mathcal C(s)=\sum c_nn^{-s}$. Then
\[
  \sum_{0<\gamma\le T}\mathcal A(\rho)\overline{\mathcal C(\rho)}=N(T)\sum_n\frac{a_nc_n}{n}-\frac{T}{2\pi}\sum_{\ell\ge2}\sum_{n\ge1}\frac{\Lambda(\ell)}{\ell n}\bigl(a_{\ell n}c_n+a_nc_{\ell n}\bigr)+O\bigl(K^2M\log^3(2T)\bigr).
\]
\end{lemma}

\begin{proof}
By~\eqref{eq:conj}, $\mathcal A(\rho)\overline{\mathcal C(\rho)}=\sum_{m,n}a_mc_nn^{-1}(n/m)^{\rho}$. The diagonal gives the first term. For $m<n$ the main term of Theorem~\ref{thm:LG} is $-\frac T{2\pi}\Lambda(n/m)$. For $m>n$ we write $(n/m)^{\rho}=(n/m)\overline{(m/n)^{\rho}}$ and obtain the main term $-\frac{T}{2\pi}\Lambda(m/n)a_mc_n/m$. Together these give the second term.

For the errors, let $s<b\le M$ be the smaller and the larger index and $X=b/s$. The term of the pair is $b^{-1}$ times at most $K^2$ times $\sum_\gamma X^\rho$ or its conjugate. The three error terms of Theorem~\ref{thm:LG}, divided by $b$, are $e_1\le2s^{-1}\log(2T)\log\log(3T)$, $e_2=\log X/\delta(s,b)$ with $\delta(s,b)=s\langle b/s\rangle$, and $e_3\le\log(2T)/(b-s)$, using $\log(b/s)\ge(b-s)/b$. Clearly $\sum_{s<b\le M}(e_1+e_3)\ll M\log^2(2T)\log\log(3T)$. For $e_2$, note that $\delta(s,b)=\min_q|b-sq|$ over prime powers $q\ne b/s$ is a positive integer, and by Bertrand's postulate the minimum is attained at some $q\le2X\le2M/s$. Hence
\[
  \sum_{s<b\le M}\frac1{\delta(s,b)}\le\sum_{s\le M}\ \sum_{q\le2M/s}\ \sum_{\substack{b\le M\\b\ne sq}}\frac1{|b-sq|}\le4M(1+\log M)^2,
\]
and $\sum e_2\ll M\log^3(2M)$.
\end{proof}

After normalisation by $N(T)\asymp TL$, the error in Lemma~\ref{lem:errorsum} is $O(K^2ML^2/T)$, which is $o(1)$ when $M\le T^{1-\varepsilon}$. It is the length of the \emph{longer} of the two polynomials that matters, not the product of the lengths. In particular, if $U_1,\dots,U_4$ are polynomials of length $T^\theta$ with bounded coefficients, $\theta<\frac12$, then $\E U_1U_2\overline{U_3U_4}$ is given by the main terms of Lemma~\ref{lem:errorsum} up to $o(1)$. By partial summation over $\gamma\in[T,2T]$, the same holds for $\EW$.

\subsection{The theorem of Benli, Elma and Ng}
Let $X(s)=\sum_{n\le N}x(n)n^{-s}$ and $Y(s)=\sum_{n\le N}y(n)n^{-s}$. Benli, Elma and Ng~\cite{BEN24} evaluate $\sum_\gamma\zeta'(\rho)X(\rho)Y(1-\rho)$ unconditionally, under the following assumptions~\cite[(4), (6)--(8)]{BEN24}:
\begin{enumerate}[label=(B\arabic*),leftmargin=*]
\item\label{BEN4} $N\ll T^{\vartheta}$ for some fixed $\vartheta<\frac12$, and $x(n)=y(n)=0$ for $n>N$;
\item\label{BEN7} $x(mn)\ll|x(m)x(n)|$ and $y(mn)\ll|y(m)y(n)|$ for all $m,n$;
\item\label{BEN8} $x(n)\ll\tau_{k_1}(n)(\log n)^{\ell_1}$ and $y(n)\ll\tau_{k_2}(n)(\log n)^{\ell_2}$ for some fixed $k_i,\ell_i$.
\end{enumerate}
Write
\[
  \mathcal F_{\alpha,h,k}=\frac{\1_{k=1}}{h^{\alpha}}\frac{\zeta'}{\zeta}(1+\alpha)-\frac{\Lambda(k)}{h^{\alpha}\prod_{p\mid k}(1-p^{-1-\alpha})}-\frac{k}{\varphi(k)}\prod_{p\mid k}(1-p^{-\alpha})\,\zeta(1-\alpha)\frac{(T/2\pi k)^{-\alpha}}{1-\alpha},
\]
which is $2\pi/T$ times~\cite[(10)]{BEN24}, and $K_{h,k}=\partial_\alpha\mathcal F_{\alpha,h,k}|_{\alpha=0}$.

\begin{theorem}[Benli--Elma--Ng {\cite[Theorem~1.1 and Corollary~1.2]{BEN24}}]\label{thm:BEN}
Under \ref{BEN4}--\ref{BEN8}, for every $A>0$,
\begin{multline*}
  \sum_{0<\gamma\le T}\zeta'(\rho)X(\rho)Y(1-\rho)=\frac{T}{2\pi}\Bigl[-(L-1)\sum_n\frac{(\log*x)(n)y(n)}{n}+\sum_n\frac{(\Lambda*\log*x)(n)y(n)}{n}\\
  +\sum_{g}\sum_{(h,k)=1}\frac{y(gh)x(gk)}{ghk}K_{h,k}\Bigr]+O_A\bigl(T(\log T)^{-A}\bigr).
\end{multline*}
\end{theorem}

This is the derivative at $\alpha=0$ of~\cite[(9)]{BEN24}, which is how~\cite[Corollary~1.2]{BEN24} is proved. The coefficient of the $\Lambda*\log*x$ term is $+T/2\pi$, as in~\cite[(26)]{BEN24}; the specialised display~\cite[(27)]{BEN24} has a misprint in this coefficient and in one Laurent coefficient, and we do not use it. Our coefficients depend on $T$, through their length and through $L$. An inspection of~\cite{BEN24} shows that the implied constants depend only on $A$, $\vartheta$, $k_i,\ell_i$ and the implied constants in~\ref{BEN4}--\ref{BEN8}: the coefficients enter the error terms only through upper bounds for $|x(n)|$ and $|y(n)|$. In~\cite{BEN24}, $\ll$ in~\ref{BEN8} refers to large $n$ only.

\begin{lemma}\label{lem:kernels}
Let $h\le T$ and let $k$ be squarefree with $(h,k)=1$. Write $\lambda=\log\ell$ and $\lambda_i=\log\ell_i$. Then
\begin{enumerate}[label=\textup{(\alph*)}]
\item $K_{h,1}=\tfrac12(L^2-\log^2h)+O(L)$;
\item $K_{h,\ell}=\frac{\ell}{\ell-1}\lambda\bigl(\log h-L+\tfrac12\lambda\bigr)+O\bigl(\lambda\bigr)$ for a prime $\ell$;
\item $K_{h,\ell_1\ell_2}=\frac{k}{\varphi(k)}\lambda_1\lambda_2$ for distinct primes $\ell_1,\ell_2$;
\item $K_{h,k}=0$ if $\omega(k)\ge3$.
\end{enumerate}
\end{lemma}

\begin{proof}
Expand $\zeta(1-\alpha)=-\alpha^{-1}+\gamma_0+O(\alpha)$, $\frac{\zeta'}{\zeta}(1+\alpha)=-\alpha^{-1}+\gamma_0+O(\alpha)$ and $1-p^{-\alpha}=\alpha\log p-\frac12\alpha^2\log^2p+O(\alpha^3)$. If $\omega(k)\ge3$, then $\Lambda(k)=0$ and the last term of $\mathcal F$ is $O(\alpha^2)$, which gives~(d). If $k=\ell_1\ell_2$, only the last term contributes, and it is $\frac k{\varphi(k)}\alpha\lambda_1\lambda_2+O(\alpha^2)$. If $k=\ell$, the second term is $-\frac{\ell\lambda}{\ell-1}+\alpha\frac{\ell\lambda\log h}{\ell-1}+O(\alpha\lambda^2/\ell)$, and the third is $\frac{\ell}{\ell-1}\bigl(\lambda-\alpha(\tfrac12\lambda^2+(\gamma_0-1)\lambda+(L-\lambda)\lambda)\bigr)+O(\alpha^2)$. Adding gives~(b). Part~(a) is the same computation for $k=1$ to second order; the poles cancel.
\end{proof}

\subsection{Prime sums}
Let $u_p=\log p/L$. By Mertens' theorem the measures $\sum_p\frac{\log p}{Lp}\delta_{u_p}$ converge weakly to Lebesgue measure on every compact interval. Hence, for every $d\ge1$ and every bounded $G$ on $[0,\infty)^d$ with compact support that is continuous outside finitely many hyperplanes,
\begin{equation}\label{eq:mertens}
  \frac1{L^d}\sum_{\substack{p_1,\dots,p_d\\\text{distinct}}}\ \prod_{i=1}^d\frac{\log p_i}{p_i}\,G(u_{p_1},\dots,u_{p_d})\ \to\ \int_{[0,\infty)^d}G .
\end{equation}
The tuples with a repeated prime contribute $O(L^{-1})$. When $G$ is Lipschitz the rate is $O(L^{-1})$.

\subsection{Character sums}
We use two classical results~\cite[Theorem~7.13 and Corollary~5.29]{IK04}.

\begin{theorem}[Multiplicative large sieve]\label{thm:LS}
For complex $a_n$,
\[
  \sum_{r\le R}\frac{r}{\varphi(r)}\ \sideset{}{^*}\sum_{\chi\bmod r}\Bigl|\sum_{n\le N}a_n\chi(n)\Bigr|^2\le(R^2+N)\sum_{n\le N}|a_n|^2,
\]
where $\sum^*$ runs over primitive characters.
\end{theorem}

\begin{theorem}[Siegel--Walfisz]\label{thm:SW}
For all $A,B>0$, uniformly for primitive $\chi\ne\chi_0$ of conductor $r\le(\log x)^B$,
\[
  \sum_{p\le x}\chi(p)\log p\ll_{A,B}x(\log x)^{-A}.
\]
The implied constant is ineffective.
\end{theorem}

\section{The local multiplicity inequality}\label{sec:local}

\subsection{Kernels and the witness}\label{sec:kernels}
Let $g$ be a real, even, piecewise polynomial function on $\mathbb R$, supported on $[-\frac12,\frac12]$, with $\int g=1$. Put $e(x)=e^{2\pi ix}$ and
\[
  k=g*g,\qquad \check g(y)=\int g(z)e(zy)\,dz,\qquad K(x)=\check g(x)^2 .
\]
Then $K\ge0$ on $\mathbb R$, $K(0)=1$, and the Fourier transform of $K$ is $k$, which is even, supported on $[-1,1]$, and a polynomial on $[0,1]$. Also $k(0)=\int g^2\ge1$. For $0<\psi<1$ we dilate:
\[
  K_\psi(x)=K(\psi x),\qquad k_\psi(u)=\psi^{-1}k(u/\psi),
\]
so that $\widehat{K_\psi}=k_\psi$ is supported on $[-\psi,\psi]$. The \emph{witness} and its \emph{offset} are
\begin{equation}\label{eq:Jdef}
  J(s)=\frac1L\sum_{n}\frac{\Lambda(n)}{n^{s}}k_\psi\Bigl(\frac{\log n}L\Bigr),\qquad a=1-\frac{k_\psi(0)}2 .
\end{equation}
$J$ is a Dirichlet polynomial of length $e^{\psi L}=(T/2\pi)^\psi$ with real coefficients, supported on prime powers.

For the Fej\'er kernel $g=\1_{[-1/2,1/2]}$ we have $k(u)=1-|u|$ and $K(x)=(\sin\pi x/\pi x)^2$. The kernels we use are even polynomials on $[-\frac12,\frac12]$, written as $g(z)=\sum_jg_j\Leg_{2j}(2z)$ with $g_0=1$. The second-order certificate of Theorem~\ref{thm:second} uses
\begin{equation}\label{eq:g8}
  g_8(z)=1-\tfrac{12}{125}\Leg_2(2z)+\tfrac{17}{1000}\Leg_4(2z),
\end{equation}
for which, on $0\le u\le1$,
\[
  k_8(u)=\tfrac{45084389}{45000000}-\tfrac{848241}{1000000}u-\tfrac{8089}{31250}u^2-\tfrac{86471}{150000}u^3+\tfrac{119}{100}u^4-\tfrac{1312253}{2500000}u^5+\tfrac{4709}{250000}u^7-\tfrac{2023}{900000}u^9 .
\]
Compared with the Fej\'er kernel, $g_8$ puts less weight near the edges of $[-\frac12,\frac12]$. This makes $K$ decay faster, at the cost of a slightly larger $k(0)$.

\subsection{The inequality}
\begin{lemma}\label{lem:local}
Assume RH, and fix $g$ and $0<\psi<1$. For every zero $\rho=\frac12+i\gamma$ with $T\le\gamma\le2T$,
\begin{equation}\label{eq:local}
  m_\rho\ \le\ \sum_{\gamma'}K_\psi\Bigl(\frac{(\gamma'-\gamma)L}{2\pi}\Bigr)=k_\psi(0)\frac{\log(\gamma/2\pi)}{L}-2\re J(\rho)+O\Bigl(\frac{\log T}{TL}\Bigr),
\end{equation}
where $\gamma'$ runs over the ordinates of all nontrivial zeros, with multiplicity, and the implied constant depends only on $g$ and $\psi$. In particular, there is a constant $c_0=c_0(g,\psi)$ such that, for $T\le\gamma\le2T$,
\begin{equation}\label{eq:halfplane}
  \re\bigl(J(\rho)+a\bigr)\le\frac{c_0}L\quad\text{if }m_\rho\ge2,\qquad
  \re\bigl(J(\rho)+a\bigr)\le\frac12+\frac{c_0}L\quad\text{if }m_\rho=1 .
\end{equation}
\end{lemma}

The idea is simple. All the ordinates are real, so every term of the sum is nonnegative, and the zero $\rho$ itself contributes $m_\rho K_\psi(0)=m_\rho$. The explicit formula computes the sum.

\begin{proof}
We use the Guinand--Weil explicit formula in the form~\cite[Lemma~1]{GG07} (see also~\cite[Theorem~5.12]{IK04}). Let $h$ be holomorphic in the strip $|\im r|\le1$, real on $\mathbb R$, with $h(r)\ll(1+|r|)^{-2}$ there, and put $\hat h(x)=\int h(r)e^{-irx}dr$. Under RH,
\begin{multline*}
  \sum_{\gamma'}h(\gamma')=h\bigl(\tfrac i2\bigr)+h\bigl(-\tfrac i2\bigr)+\frac1{2\pi}\int_{\mathbb R}h(r)\Bigl(\re\frac{\Gamma'}{\Gamma}\Bigl(\frac14+\frac{ir}2\Bigr)-\log\pi\Bigr)dr\\-\frac1{2\pi}\sum_{n}\frac{\Lambda(n)}{\sqrt n}\bigl(\hat h(\log n)+\hat h(-\log n)\bigr).
\end{multline*}
Take $h(r)=\check g(\psi z)\check g(-\psi z)$ with $z=(r-\gamma)L/2\pi$. This is entire, and equals $K_\psi((r-\gamma)L/2\pi)$ on $\mathbb R$ because $g$ is even. Integrating by parts once, $|\check g(y)|\ll e^{\pi|\im y|}(1+|y|)^{-1}$, so $h$ is admissible and $|h(\pm\frac i2)|\ll T^{\psi/2-2}$. Next, $\hat h(x)=\frac{2\pi}Le^{-i\gamma x}k_\psi(x/L)$, so the prime sum equals $-2\re J(\rho)$ exactly. In the gamma integral substitute $r=\gamma+2\pi x/L$. By Stirling's formula $\re\frac{\Gamma'}{\Gamma}(\frac14+\frac{iv}2)=\log\frac{|v|}2+O(v^{-2})$ for $|v|\ge1$, and it is $O(\log(2+|v|))$ for all real $v$. Since $K_\psi(x)\ll(1+|x|)^{-2}$ and $\int K_\psi=k_\psi(0)$, the gamma term is $k_\psi(0)\log(\gamma/2\pi)/L+O(\log T/(TL))$. This proves the identity in~\eqref{eq:local}, and the inequality holds because every term is nonnegative.

For~\eqref{eq:halfplane}, note that $\log(\gamma/2\pi)/L=1+O(L^{-1})$. If $m_\rho\ge2$, then~\eqref{eq:local} gives $2\re J(\rho)\le k_\psi(0)-2+O(L^{-1})=-2a+O(L^{-1})$. If $m_\rho=1$, it gives $2\re J(\rho)\le k_\psi(0)-1+O(L^{-1})$.
\end{proof}

\begin{remark}\label{rem:localclassical}
The inequality $m_\rho\le\sum_{\gamma'}K((\gamma'-\gamma)L/2\pi)$ is classical. Montgomery~\cite{Mon73} sums it over $\gamma$ and evaluates the double sum by his pair correlation theorem; this is the basis of the pair correlation method~\cite{CG93,CGdL20}. Goldston and Gonek~\cite{GG07} use it at a single zero to bound $m_\rho$ by $\log T/\log\log T$. We use it differently: as a constraint on the \emph{location} of $J(\rho)$ at each multiple zero, which is then combined, at the same zero, with the vanishing of the detector.
\end{remark}

The constraint is informative only if $a>0$, that is, $k_\psi(0)=k(0)/\psi<2$. Since $k(0)\ge1$, this forces $\psi>\frac12$. So the witness must be longer than $T^{1/2}$, and the mixed moment with the detector involves a polynomial longer than any to which the theorem of Benli, Elma and Ng applies. This is the reason for Section~\ref{sec:longmixed}.

\section{The second-order certificate}\label{sec:second}

\subsection{The certificate inequality}
\begin{lemma}\label{lem:cert}
Assume RH. Let $F$ be a function on the zeros that vanishes at every multiple zero, and let $\tau\ge0$ be fixed. For $L>c_0\tau$,
\[
  1-\EW\1_{\mathrm s}\ \le\ \Bigl(1-\frac{c_0\tau}L\Bigr)^{-2}\,\EW\bigl|1-F-\tau(J+a)\bigr|^2 .
\]
\end{lemma}

\begin{proof}
At a multiple zero with $T\le\gamma\le2T$, $F=0$ and, by~\eqref{eq:halfplane}, $\re(1-F-\tau(J+a))\ge1-c_0\tau/L$. So $|1-F-\tau(J+a)|^2\ge(1-c_0\tau/L)^2(1-\1_{\mathrm s})$ at every zero in the support of $W(\gamma/T)$.
\end{proof}

Nothing is assumed about the values of $F$ and $J$ at the simple zeros, and no moment of order higher than two is needed.

\subsection{Optimisation}
Let $X_1,\dots,X_{n-1}$ be functions on the zeros vanishing at the multiple zeros, and let $X_n=J+a$. Suppose that, as $T\to\infty$,
\begin{equation}\label{eq:Gm}
  \EW X_i\to m_i\in\mathbb R,\qquad \re\EW X_i\overline{X_j}\to G_{ij},
\end{equation}
with $G$ positive definite. For real $\beta\in\mathbb R^n$, $\EW|1-\beta^{\mathsf T}X|^2\to1-2\beta^{\mathsf T}m+\beta^{\mathsf T}G\beta$, which is minimal at $\beta=G^{-1}m$, with value $1-m^{\mathsf T}G^{-1}m$.

\begin{proposition}\label{prop:optim}
Assume RH and~\eqref{eq:Gm}. Write $m=(A,h)$ with $A\in\mathbb R^{n-1}$, and write the limiting covariance $\Sigma=G-mm^{\mathsf T}$ in blocks $\begin{psmallmatrix}V&C\\C^{\mathsf T}&V_J\end{psmallmatrix}$, with $\Sigma$ positive definite. Put
\[
  \Ksq=A^{\mathsf T}V^{-1}A+\frac{\gap_+^2}{D_J},\qquad \gap=h-A^{\mathsf T}V^{-1}C,\qquad D_J=V_J-C^{\mathsf T}V^{-1}C .
\]
Then $\liminf_{T\to\infty}\EW\1_{\mathrm s}\ge\Ksq/(1+\Ksq)$.
\end{proposition}

\begin{proof}
By the Sherman--Morrison formula, $G^{-1}m=\Sigma^{-1}m/(1+m^{\mathsf T}\Sigma^{-1}m)$, and block elimination gives $m^{\mathsf T}\Sigma^{-1}m=A^{\mathsf T}V^{-1}A+\gap^2/D_J$. The last entry of $G^{-1}m$ is $\tau=\gap/(D_J(1+m^{\mathsf T}\Sigma^{-1}m))$. If $\gap\ge0$, apply Lemma~\ref{lem:cert} with $F=\sum_{i<n}\beta_iX_i$ and this $\tau$. If $\gap<0$, drop $X_n$; the bound is then $A^{\mathsf T}V^{-1}A/(1+A^{\mathsf T}V^{-1}A)$.
\end{proof}

The proposition says how the two tests combine. The detectors alone give $\Ksq_0/(1+\Ksq_0)$ with $\Ksq_0=A^{\mathsf T}V^{-1}A$; for one detector this is the quotient $A^2/(A^2+V)$ of Bui and Heath-Brown. The witness adds $\gap^2/D_J$. Here $A^{\mathsf T}V^{-1}C$ is the part of the mean of $J+a$ that the detectors can predict through their covariance with $J$, and $\gap$ is the part they cannot. So the witness helps exactly when the mean of $J+a$ exceeds its prediction from the detectors. The denominator $D_J$ is the variance of $J$ left after projecting onto the detectors.

\subsection{Moments of the witness}
\begin{lemma}\label{lem:Jmom}
Assume RH, and let $0<\psi<1$. Let $J_{\mathrm P}$ be the part of $J$ supported on primes. Then
\[
  \EW J\to\mu:=-\int_0^\psi uk_\psi(u)\,du,\qquad \EW|J|^2\to\int_0^\psi uk_\psi(u)^2du,\qquad \EW|J-J_{\mathrm P}|^2=O(L^{-2}).
\]
\end{lemma}

\begin{proof}
The coefficients of $J$ are bounded and supported on $n\le T^\psi$, so Lemma~\ref{lem:errorsum} applies. With $c_n=\1_{n=1}$, only $n=1$ and $\ell$ a prime power occur, and
\[
  \EW J=-\frac1{L^2}\sum_\ell\frac{\Lambda(\ell)^2}{\ell}k_\psi\Bigl(\frac{\log\ell}L\Bigr)+o(1)\to-\int_0^\psi uk_\psi(u)\,du .
\]
With $c=a$, the diagonal gives $L^{-2}\sum_n\Lambda(n)^2k_\psi(\log n/L)^2/n\to\int uk_\psi^2$. Off the diagonal, $a_{\ell n}a_n\ne0$ forces $n=p^j$, $\ell n=p^i$ with $i>j\ge1$, and these terms are $O(L^{-3})$. The same computation for $J-J_{\mathrm P}$ gives $O(L^{-2})$.
\end{proof}

The prime powers cost nothing in the pointwise inequality~\eqref{eq:local}, and in averages they can be dropped: by Cauchy--Schwarz and Lemma~\ref{lem:BHBW}, $\EW F_P\overline{(J-J_{\mathrm P})}=O(L^{-1})$.

\subsection{The mixed moment}
The remaining entries of $G$ are the mixed moments $\re\EW F_P\overline J$. They are given by the case $d=1$ of Theorem~\ref{thm:longmixed} in the next section. With Lemma~\ref{lem:Jmom} and the centring $\EW F_P\,\EW\bar J\to A_P\mu$, it gives:

\begin{corollary}\label{cor:C}
Assume RH. Let $0<\theta<\frac12$, $0<\psi<1$, and $P(0)=0$, $f=f_P$. Then
\[
  \re\EW\bigl[(F_P-\EW F_P)\overline{(J-\EW J)}\bigr]\to C_P:=\int_0^\psi k_\psi(u)\,\mathfrak c_f(u)\,du,
\]
where
\[
  \mathfrak c_f(u)=\frac{f(0)-f(u)}2+\int_0^uf-\frac{f(0)}2(u+u^2).
\]
\end{corollary}

For the Bui--Heath-Brown taper at the endpoint this is the function $\mathfrak c$ of Section~\ref{sec:sketch}. Note that $\mathfrak c_f(u)$ is a polynomial in $u$ for $u\le\theta$, and of the form $c_0+c_1u+c_2u^2$ for $u\ge\theta$; the long primes $u>\theta$ contribute through $f(0)$ and $I_f$ only.

\subsection{Proof of the second-order bound}
Fix $\theta<\frac12<\psi<1$, a kernel $g$, and polynomials $P_0,\dots,P_{n-2}$ vanishing at $0$. Take $X_{i+1}=F_{P_i}$ and $X_n=J+a$. By Lemmas~\ref{lem:BHBW} and~\ref{lem:Jmom} and Corollary~\ref{cor:C}, the limits~\eqref{eq:Gm} exist, with
\[
  m=\begin{pmatrix}(A_{P_i})_i\\ \mu+a\end{pmatrix},\qquad \Sigma=\begin{pmatrix}(V_{P_iP_j})_{ij}&(C_{P_i})_i\\ (C_{P_j})_j^{\mathsf T}&\int_0^\psi uk_\psi^2-\mu^2\end{pmatrix}.
\]
All entries are integrals of piecewise polynomials. They are rational when $\theta,\psi$ are rational, and continuous in $(\theta,\psi)$ up to the endpoint $(\frac12,1)$. For a fixed real vector $\beta$ with $\beta_n\ge0$, put $\kappa(\theta,\psi;\beta)=2\beta^{\mathsf T}m-\beta^{\mathsf T}G\beta$. Lemma~\ref{lem:cert} and Lemma~\ref{lem:unweight} give $\liminf N_s/N\ge\kappa(\theta,\psi;\beta)$ for every admissible $(\theta,\psi)$, and letting $(\theta,\psi)\to(\frac12,1)$,
\begin{equation}\label{eq:supcert}
  \liminf_{T\to\infty}\frac{N_s(T)}{N(T)}\ \ge\ \kappa\bigl(\tfrac12,1;\beta\bigr).
\end{equation}
The limit $T\to\infty$ is always taken first, at fixed exponents. No estimate is needed uniformly as $\theta\to\frac12$ or $\psi\to1$.

\begin{theorem}\label{thm:second}
Assume RH. Then $\liminf_{T\to\infty}N_s(T)/N(T)\ge0.7073645768999$.
\end{theorem}

\begin{proof}
Take the kernel $g_8$ of~\eqref{eq:g8} and the twelve tapers $b_0,\dots,b_{11}$ of~\eqref{eq:basis}. At $(\theta,\psi)=(\frac12,1)$, the vector $A$ and the matrix $V$ are given by~\eqref{eq:basisAV}, and
\[
  \mu=-\tfrac{45300019}{277200000},\qquad a=\tfrac{44915611}{90000000},\qquad V_J=\tfrac{3535240493473142893}{62432370000000000000}.
\]
The twelve covariances $C_{P_j}$ are rational numbers, computed from Corollary~\ref{cor:C}; they are listed in Appendix~\ref{app:data}. Then
\[
  \gap=\tfrac{352643608931}{21086208000000}>0,\qquad
  \tau=\tfrac{4975789992282542264832000000000}{6735054433292025624721100070697}>0,
\]
and Proposition~\ref{prop:optim} with~\eqref{eq:supcert} gives
\[
  \liminf_{T\to\infty}\frac{N_s(T)}{N(T)}\ \ge\ \frac{\Ksq}{1+\Ksq}=\frac{42877250366431111605585433888831}{60615489899628230622489900636273}=0.70736457689990\ldots .\qedhere
\]
\end{proof}

\subsection{Remarks on the certificate}
\emph{The endpoint is not essential.} At the fixed exponents $(\theta,\psi)=(\frac{4999}{10000},\frac{9999}{10000})$ the coefficients $\beta$ of the proof give $0.70730464947\ldots$, already above $19/27$ by $0.0036$.

\emph{The full taper space.} The tapers $b_j$ span all polynomials vanishing at $\frac12$. For a fixed kernel, the optimum of Proposition~\ref{prop:optim} over all tapers is a finite-dimensional problem: the functionals $A$ and $C$ are represented, for the inner product $V$, by polynomials, so the optimum is attained by a polynomial taper of bounded degree. For the kernel
\begin{equation}\label{eq:gstar}
  g_*(z)=1+10^{-6}\bigl(-95913\Leg_2(2z)+16977\Leg_4(2z)-1307\Leg_6(2z)+244\Leg_8(2z)\bigr)
\end{equation}
the optimum is
\begin{equation}\label{eq:p0}
  p_0=\frac{3413914900541438580624509235808381477124808899103359}{4826245010227457178772530978845119244788120650314257}=0.70736460608752\ldots,
\end{equation}
with witness coefficient $\tau_*=0.73879626316\ldots$. This certificate is used in Section~\ref{sec:distinct}. Proposition~\ref{prop:ceiling} below shows that no kernel and taper can do much better.

\emph{Why $\psi>\theta$ matters.} If the witness is no longer than the mollifier, it only sees primes that the detector already sees, and the gap is negative. For example, for the Fej\'er kernel with $\psi=\theta=\frac12$ we have $a=0$, $\mu=-\frac1{12}$, and with the taper of Section~\ref{sec:sketch}, $C=\frac{19}{240}$ and $\gap=-\frac{77}{240}$. Then $\tau=0$ and the witness gives nothing. The gain comes from the primes in $(T^{1/2},T)$.

\section{The mixed moment with a long numerator}\label{sec:longmixed}

\subsection{Statement}
For $n=p_1\cdots p_d$ squarefree write $u(n)=u_{p_1}\cdots u_{p_d}$, where $u_p=\log p/L$. For $S\subseteq\{1,\dots,d\}$ and $\mathbf u\in[0,\infty)^d$ put
\[
  u_S=\sum_{i\in S}u_i,\qquad u_h=\sum_{i\notin S}u_i,\qquad V_S=\sum_{i\notin S}u_i^2,
\]
and define the density
\begin{equation}\label{eq:Dd}
  D_d(\mathbf u)=\sum_{S\subseteq\{1,\dots,d\}}(-1)^{|S|}\Bigl[f(u_S)\Bigl(\frac12-u_h+\frac{u_h^2}2-V_S\Bigr)+(1-u_h)\,I_f(u_S)\Bigr].
\end{equation}
Here $f=f_P$ is extended by $0$ beyond $\theta$, and $I_f(x)=\int_x^\theta f$ is $0$ for $x\ge\theta$. The letter $h$ in $u_h$ stands for the part of the numerator that is not shared with the mollifier variable; this will be clear from the proof.

\begin{theorem}\label{thm:longmixed}
Assume RH. Let $0<\theta<\frac12$, $0<\psi<1$ and $d\ge1$. Let $\mathcal K$ be a bounded symmetric function on $[0,\infty)^d$, continuous outside finitely many hyperplanes and supported in $\{u_1+\dots+u_d\le\psi\}$, and put
\[
  \mathcal A(s)=\sum_{p_1<\dots<p_d}u_{p_1}\cdots u_{p_d}\,\mathcal K(u_{p_1},\dots,u_{p_d})\,(p_1\cdots p_d)^{-s}.
\]
Let $P$ be a real polynomial with $P(0)=0$, and $f=f_P$. Then
\[
  \EW F_P\overline{\mathcal A}\ \to\ \frac1{d!}\int_{[0,\infty)^d}\mathcal K(\mathbf u)\,D_d(\mathbf u)\,d\mathbf u .
\]
\end{theorem}

For $d=1$ the density is
\begin{equation}\label{eq:D1}
  D_1(u)=f(0)\Bigl(\frac12-u-\frac{u^2}2\Bigr)+(1-u)I_f-\frac{f(u)}2-\int_u^\theta f .
\end{equation}
The prime part of the witness is $J_{\mathrm P}=\mathcal A$ with $d=1$ and $\mathcal K=k_\psi$. Corollary~\ref{cor:C} follows from~\eqref{eq:D1}, since $\EW F_P\,\EW\bar J\to A_P\mu$ and $D_1(u)+uA_P=\mathfrak c_f(u)$.

When $\psi\le\theta$ the polynomial $\mathcal A$ is short, and Theorem~\ref{thm:longmixed} is a consequence of the theorem of Benli, Elma and Ng; this is Step~2 below. The new case is $\psi>\theta$, and in particular $\psi>\frac12$.

\subsection{Why the long numerator is harmless}
Consider $d=1$, so that $\mathcal A$ is a sum over primes $p\le T^\psi$. Following Conrey, Ghosh and Gonek, one writes the sum over zeros as a contour integral and moves to the left of the critical strip with the functional equation. Gonek's lemma then turns the integral into a sum of additive characters,
\[
  \sum_{h\le T^\theta}\frac{b(h)}h\sum_{p}y_p\sum_{m\asymp hT/p}e\Bigl(\frac{pm}h\Bigr)\cdot(\text{smooth weight}),
\]
where $b(h)$ are the mollifier coefficients and $y_p$ those of $\mathcal A$. The denominator of the additive character is $h$, the mollifier variable. It is at most $T^\theta$ whatever the length of $\mathcal A$. The long prime $p$ sits in the numerator. Expanding in Dirichlet characters modulo divisors of $h$, the principal characters give the main term. The others lead to character sums $\sum_p y_p\chi(p)p^{-it}$ over primes, of conductor at most $T^\theta$. These have cancellation on average over $\chi$ by the large sieve, and individually, for small conductors, by Siegel--Walfisz.

Compare this with the barrier $\theta<\frac12$ for the mollifier itself. In the second moment $\sum|\zeta'B|^2$ the mollifier occurs twice; after the same steps the modulus is a mollifier variable and the frequency is the other one. A mollifier of length $T^{1/2+\delta}$ then produces moduli beyond $T^{1/2}$, and the large sieve term $R\sqrt T$ in~\eqref{eq:LSblock} below exceeds $T$. In the mixed moment only the mollifier becomes a modulus, and the length $T^\psi$ of $\mathcal A$ enters only through the terms $\sqrt{TP_0}$ and $T\sqrt{R/P_0}$, which are harmless for $\psi<1$.

There is one complication. A long numerator $n=p_1\cdots p_d$ may share primes with $h$, and after removing the shared primes what remains may be small. Then there is no long prime left to sum over, and Siegel--Walfisz gives nothing. We therefore cut $\mathcal A$ by the size of $n$ at $T^\beta$, with $\theta<\beta<\frac12$. Below the cut we use the theorem of Benli, Elma and Ng. Above it, the unshared part of $n$ exceeds $T^{\beta-\theta}$, so it has a prime factor at least $T^{(\beta-\theta)/d}$.

\subsection{Proof of Theorem~\ref{thm:longmixed}}
Write $Y=T^\theta$, $X_0=T/2\pi$, $b(h)=\mu(h)f(\log h/L)$ and $y_n=u(n)\mathcal K(\mathbf u_n)$ for squarefree $n$ with $d$ prime factors, where $\mathbf u_n$ is the vector of the $u_p$, $p\mid n$. By~\eqref{eq:conj} we must evaluate
\begin{equation}\label{eq:Sigma}
  \Sigma:=\frac1{N_WL}\sum_{\gamma>0}W\Bigl(\frac\gamma T\Bigr)\zeta'(\rho)B_P(\rho)\mathcal A(1-\rho).
\end{equation}
All error terms are for fixed $\theta,\psi,\beta,d,P,\mathcal K,W$.

\emph{Step 0: reductions.} Both sides of the theorem are linear in $\mathcal K$. Each side changes by $O(\|\mathcal K-\mathcal K'\|)$, with the norm of $L^2(\prod_iu_i\,du_i)$, if $\mathcal K$ is replaced by $\mathcal K'$. For the right side this is clear since $D_d$ is bounded. For the left side use Cauchy--Schwarz, Lemma~\ref{lem:BHBW}, and the fact that $\EW|\mathcal A|^2\to\frac1{d!}\int\mathcal K^2\prod u_i\,du_i$, by Lemma~\ref{lem:errorsum} and~\eqref{eq:mertens}. (For a homogeneous polynomial of degree $d$ over squarefree numbers, the off-diagonal terms of Lemma~\ref{lem:errorsum} vanish, since $\ell n$ and $n$ cannot both have $d$ prime factors.) Hence we may assume that $\mathcal K$ is a finite linear combination of indicator functions of boxes, symmetrised, with support in $\{u_1+\dots+u_d\le\psi'\}$ for some $\psi'<1$. We write $\psi$ for $\psi'$.

\emph{Step 1: the cut.} Fix $\beta$ with $\theta<\beta<\frac12$, and write $\mathcal A=\mathcal A_\le+\mathcal A_>$, where $\mathcal A_\le$ contains the terms with $u_1+\dots+u_d\le\beta$. Both parts are of the same form as $\mathcal A$, with $\mathcal K$ multiplied by the indicator of a half-space. We prove the theorem for each of them.

\emph{Step 2: short numerators.} We apply Theorem~\ref{thm:BEN} with $N=T^\beta$, $x(h)=b(h)$ and $y(n)$ the coefficients of $\mathcal A_\le$. Conditions~\ref{BEN4} and~\ref{BEN8} hold. Condition~\ref{BEN7} need not hold, since $x$ and $y$ may change sign and $y(1)=0$. We use a background. Let
\[
  B_0(n)=\1\bigl[n\le T^\beta,\ \mu^2(n)=1,\ \omega(n)\le d\bigr]\,u(n),\qquad f_0(v)=(1-v/\theta)_+ .
\]
The support of $B_0$ is closed under taking divisors, and $B_0(mn)=B_0(m)B_0(n)$ when $mn$ lies in it; so $B_0$ satisfies~\ref{BEN7}. Since $|y|\le CB_0$ for some $C$, so does $y+(C+1)B_0$, which is at least $B_0$. Similarly $f_0(s+t)\le f_0(s)f_0(t)$, and $|f|\le C'f_0$ because $f(\theta)=P(0)=0$; so $\mu(h)f_0(\log h/L)$ and $\mu(h)(f+(C'+1)f_0)(\log h/L)$ satisfy~\ref{BEN7}. The main terms of Theorem~\ref{thm:BEN} are bilinear in $(x,y)$, so four applications give the main terms for $(x,y)$, with error $O(T(\log T)^{-A})$. As noted after Theorem~\ref{thm:BEN}, the error terms are uniform in our $T$-dependent coefficients. Partial summation over $\gamma\in[T,2T]$ introduces the weight $W$.

It remains to evaluate the main terms. Let $n=p_1\cdots p_d$ be in the support of $y$, and write $g=\prod_{i\in S}p_i$ for a divisor of $n$, so that $x(g)=(-1)^{|S|}f(u_S)$ and $\log(n/g)=Lu_h$. For squarefree $m$ we have $(\Lambda*\log)(m)=\log^2m-\sum_{p\mid m}\log^2p$. Hence
\[
  -(L-1)(\log*x)(n)=-L^2\sum_S(-1)^{|S|}f(u_S)u_h+O(L),
\]
\[
  (\Lambda*\log*x)(n)=L^2\sum_S(-1)^{|S|}f(u_S)(u_h^2-V_S).
\]
In the last sum of Theorem~\ref{thm:BEN} write $n=gh$ with $g=\prod_{i\in S}p_i$; then $k$ runs over squarefree integers coprime to $n$, and $x(gk)=(-1)^{|S|}\mu(k)f(u_S+\log k/L)$. By Lemma~\ref{lem:kernels} only $\omega(k)\le2$ contributes. With $v=\log\ell/L$, the terms $k=1$, $k=\ell$ and $k=\ell_1\ell_2$ give, after division by $L^2$ and summation over $\ell$ by~\eqref{eq:mertens},
\[
  (-1)^{|S|}\Bigl[\frac{f(u_S)}2(1-u_h^2)+\int_0^\infty\Bigl(1-u_h-\frac v2\Bigr)f(u_S+v)\,dv+\frac12\int_0^\infty vf(u_S+v)\,dv\Bigr]+O(L^{-1}).
\]
The at most $d$ primes $\ell\mid n$ excluded from the sums over $\ell$ cost $O(L^{-1})$. The sum of the three contributions is the $S$-summand of~\eqref{eq:Dd}, since $\int_0^\infty f(u_S+v)\,dv=I_f(u_S)$. Finally $y(n)/n=\prod_i(\log p_i/Lp_i)\mathcal K(\mathbf u_n)$, and~\eqref{eq:mertens} turns the sum over $p_1<\dots<p_d$ into $\frac1{d!}\int\mathcal K D_d$.

From now on $\mathcal A=\mathcal A_>$, so every $n$ in its support has $T^\beta<n\le T^\psi$.

\emph{Step 3: the contour.} To use an entire weight, put $W_T(z)=\pi^{-1/2}\int_{\mathbb R}W(v/T)e^{-(v-z)^2}dv$. On $\mathbb R$, $W_T(t)=W(t/T)+O(T^{-2})$; in every fixed horizontal strip $W_T^{(j)}(z)\ll_jT^{-j}$, and $W_T$ decays like a Gaussian outside $[\frac T2,3T]$. Replacing $W(\gamma/T)$ by $W_T(\gamma)$ in~\eqref{eq:Sigma} costs $O(T^{-1+\varepsilon})$. The function
\[
  \frac{\zeta'}\zeta(s)\,\zeta'(s)\,B_P(s)\,\mathcal A(1-s)\,W_T\bigl((s-\tfrac12)/i\bigr)
\]
has residue $\zeta'(\rho)B_P(\rho)\mathcal A(1-\rho)W_T(\gamma)$ at a simple zero $\rho$, and is regular at multiple zeros, where $\zeta'$ vanishes. Integrate it around the rectangle with vertical sides $\re s=c:=1+1/L$ and $\re s=1-c$. The horizontal sides tend to $\pm i\infty$, and the pole at $s=1$ carries the factor $W_T(-i/2)$, which is exponentially small. So $N_WL\,\Sigma$ is the right vertical integral minus the left one, up to $O(T^{-1+\varepsilon})$.

\emph{Step 4: the right line.} On $\re s=c$, expand $(\zeta'/\zeta)\zeta'B_P(s)=\sum_m(\Lambda*\log*b)(m)m^{-s}$ and $\mathcal A(1-s)=\sum_ny_nn^{s-1}$. The term $(m,n)$ integrates to $\int W_T(t)(n/m)^{it}dt$ times bounded factors. If $m\ne n$ then $|\log(m/n)|\gg T^{-\psi}$, and repeated integration by parts makes the total of these terms $O(T^{-A})$. The diagonal gives $\frac{T}{2\pi}\int W\sum_n y_n(\Lambda*\log*b)(n)/n$. For squarefree $n$, $(\Lambda*\log*b)(n)=\sum_{g\mid n}b(g)(\Lambda*\log)(n/g)$. As in Step~2, after division by $N_WL$ and summation over the primes, this gives the density
\begin{equation}\label{eq:rightdens}
  \sum_S(-1)^{|S|}f(u_S)\bigl(u_h^2-V_S\bigr).
\end{equation}
For $d=1$ it vanishes.

\emph{Step 5: the left line.} Write $\zeta(s)=\chi(s)\zeta(1-s)$, $z_0=1-s$ and $H=\chi'/\chi(s)$. Then
\[
  \frac{\zeta'}{\zeta}(s)\zeta'(s)=\chi(s)\sum_m\frac{H^2+2H\log m+(\Lambda*\log)(m)}{m^{z_0}} .
\]
Insert $B_P(s)\mathcal A(1-s)$ and integrate term by term. The term $(h,n,m)$ has phase $t-t\log(t/2\pi r)+\pi/4$ with $r=nm/h$, and a saddle point at $t=2\pi r$, where $H=-\log r+O(T^{-1})$. By smoothed stationary phase (Gonek's lemma~\cite[Lemma~5]{Gon84}, with a smooth weight in place of the sharp cut-off) the term is $e(nm/h)$ times $W(r/X_0)$ times the bracket at the saddle, with relative error $O(T^{-1}L^C)$ when $r\asymp T$. The terms with $r$ outside $[\frac T4,4T]$ are $O(T^{-A})$. Hence the left vertical integral is $\mathcal L+O(YT^\varepsilon)$, where
\begin{equation}\label{eq:left}
  \mathcal L=\sum_{h\le Y}\frac{b(h)}h\sum_ny_n\sum_{m\ge1}e\Bigl(\frac{nm}h\Bigr)W\Bigl(\frac r{X_0}\Bigr)\Bigl[\log^2r-2\log r\log m+(\Lambda*\log)(m)\Bigr],\qquad r=\frac{nm}h .
\end{equation}
The error is $O(T^{-1+\varepsilon}\sum_{h\le Y}h^{-1}\sum_n|y_n|\,hT/n)=O(YT^\varepsilon)$. A sharp cut-off would give an error of size $T^{1/2}$ times the length of $\mathcal A$, which is too large for $\psi>\frac12$. (The sign of the phase depends on the orientation and plays no role.)

\emph{Step 6: shared primes and characters.} Fix $n$ and a squarefree $h$. Let $S$ be the set of $i$ with $p_i\mid h$, and write $g=(n,h)=\prod_{i\in S}p_i$, $h=gH$ and $n=gP$. Then $(P,gH)=1$, $e(nm/h)=e(Pm/H)$, and $b(h)=(-1)^{|S|}\mu(H)f(u_S+\log H/L)$. Let $\delta=(m,H)$, $m=\delta m'$ and $q=H/\delta$, so that $(Pm',q)=1$. Expand
\[
  e\Bigl(\frac{Pm'}q\Bigr)=\frac1{\varphi(q)}\sum_{\chi\bmod q}\tau(\bar\chi)\chi(Pm').
\]
If $\chi$ is induced by the primitive character $\chi^*$ modulo $r\mid q$, write $q=ar$; then $(a,r)=1$ and $\tau(\bar\chi)=\mu(a)\bar\chi^*(a)\tau(\bar\chi^*)$, with $|\tau(\bar\chi^*)|=\sqrt r$. The condition $(m',a)=1$ is removed by M\"obius inversion, $m'=ej$ with $e\mid a$. This rewrites $\mathcal L$ as a sum over $(g,P,\delta,a,e,r,\chi^*,j)$, bijectively. The terms with $r=1$ form the \emph{principal part}; they replace $e(Pm/H)$ by $c_H(m)/\varphi(H)$, where $c_H$ is Ramanujan's sum.

\emph{Step 7: the principal part.} For $l$ prime put $E_l(s)=(l^{1-s}-1)/(l-1)$, $E_H=\prod_{l\mid H}E_l$ and $G_H=\zeta E_H$. The Dirichlet series of $c_H(m)/\varphi(H)$ times $1$, $\log m$ and $(\Lambda*\log)(m)$ are $G_H$, $-G_H'$ and
\[
  D_H=G_H''-G_H\Bigl(\frac{\zeta'}{\zeta}\Bigr)'-\zeta(s)\sum_{l\mid H}E_{H/l}(1-E_l)\frac{l^{-s}\log^2l}{(1-l^{-s})^2}.
\]
(The $\log^2m$ part of $(\Lambda*\log)(m)=\log^2m-\sum_{l^\alpha\mid m}\alpha\log^2l$ gives $G_H''$; in the second part, conditioning on $l^\alpha\mid m$ replaces $E_l$ by $1$.) Write the $m$-sum in~\eqref{eq:left} as a Mellin integral and move it to $\re s=\frac12+\varepsilon$, keeping the residue at $s=1$. On that line the three series are $O(H^{-1/2}T^\varepsilon(1+|s|)^{O(1)})$ by RH. Since the saddle length is $X=HX_0/P$, the remainders total
\[
  \ll T^{\varepsilon}\sum_{n\le T^\psi}\frac{|y_n|}{n^{1/2}}\sum_{h\le Y}\frac{X_0^{1/2}}{h^{1/2}H^{1/2}}\ll T^{(1+\psi)/2+3\varepsilon}=o(T).
\]
For the residues, put $x=s-1$, $\lambda=\log l$, $v=\log l/L$, and note $\log X=(1-u_h+\log H/L)L+O(1)$. Normalise by $N_WL$.
\begin{itemize}[leftmargin=*]
\item $H=1$: the three Dirichlet series are $\zeta$, $-\zeta'$ and $(\zeta'/\zeta)\zeta'$, and the means of the bracket give, with the minus sign of the left line, the density $(-1)^{|S|}f(u_S)(\frac12-u_h-\frac{u_h^2}2)$.
\item $H=l$: $G_l$ and $G_l'$ are regular at $s=1$ since $E_l(1)=0$, and the polar part of $D_l$ is $\frac1l(\lambda x^{-2}-\frac32\lambda^2x^{-1})$ to leading order. With $b(gl)=-(-1)^{|S|}f(u_S+v)$ this gives the density $(-1)^{|S|}\int_0^\infty(1-u_h-\frac v2)f(u_S+v)\,dv$.
\item $H=l_1l_2$: the only polar part is $-\lambda_1\lambda_2/(\varphi(H)x)$, and the density is $(-1)^{|S|}\frac12\int_0^\infty vf(u_S+v)\,dv$.
\item $\omega(H)\ge3$: all three series are regular at $s=1$.
\end{itemize}
Lower-order Laurent coefficients, derivatives of the Mellin weight and the error in $\log X$ lose a factor $L$, and the conditions $l\nmid n$ and $l_1\ne l_2$ cost $O(L^{-1})$. Adding~\eqref{eq:rightdens}, the total density is exactly the $S$-summand of~\eqref{eq:Dd}, as in Step~2.

\emph{Step 8: large conductors.} It remains to show that the non-principal part of $\mathcal L$ is $o(TL^2)$. Fix $g,\delta,a,e$ and decompose $P\sim P_0$, $r\sim R$ and $j\sim M$ dyadically. The saddle condition $Pej/(ar)\asymp T$ forces $M\asymp aRT/(eP_0)$. Separate the variables in $W(Pej/(arX_0))$ and in the powers of $\log r$ by a Mellin integral on $\re s=0$; its transform decays rapidly, at a cost $O(L^C)$. The factors $b(g\delta ar)$, $\bar\chi^*(a)$, $\chi^*(e)$ and $(ar)^s$ depend on $r$ but are bounded, and we take absolute values in them. The prefactor $1/h$ and the Gauss sum give the weight $\ll(g\delta a^2R^{3/2})^{-1}$, up to $\log L$. What remains is
\[
  \sum_{r\sim R}\ \sideset{}{^*}\sum_{\chi^*\bmod r}\Bigl|\sum_{P\sim P_0}y_{gP}\1_{(P,g\delta a)=1}\chi^*(P)P^{-s}\Bigr|\,\Bigl|\sum_{j\sim M}c_j\chi^*(j)j^{-s}\Bigr|,
\]
with $|c_j|\ll L^2$ and $c_j$ independent of $r$. For $R\ge L^B$, Cauchy--Schwarz and Theorem~\ref{thm:LS} bound the block by
\begin{equation}\label{eq:LSblock}
  \ll L^C\Bigl(\frac{T}{g\delta ae\sqrt R}+\frac{\sqrt{TP_0}}{g\delta a^{3/2}\sqrt e}+\frac{T\sqrt{R/P_0}}{g\delta ae}+\frac{R\sqrt T}{g\delta a^{3/2}\sqrt e}\Bigr).
\end{equation}
Here $R\le Y/(g\delta a)$, $T^\beta/g<P_0\le T^\psi$ and so $R/P_0\le T^{\theta-\beta}$. The last three terms are $\ll T^{1-c}$ with $c=\min(\frac{1-\psi}2,\frac{\beta-\theta}2,\frac12-\theta)>0$, and the first is $\ll TL^{C-B/2}$. The sums over $g,\delta,a,e$ and the dyadic ranges cost $L^{O(1)}$.

\emph{Step 9: small conductors.} Let $1<R<L^B$. Since $P>T^\beta/g\ge T^{\beta-\theta}$ and $\omega(P)\le d$, the largest prime factor $p^*$ of $P$ satisfies $p^*\ge T^{(\beta-\theta)/d}$, and $p^*\nmid g\delta a$. Fix the other prime factors of $P$. Then $p^*$ runs over an interval, minus the $O(L)$ primes dividing $g\delta a$, with a weight (the factor $u_{p^*}$, the Mellin factor and $(p^*)^{-s}$) of total variation $O(L^2(1+|s|))$; here we use that $\mathcal K$ is a combination of box indicators, cut by the half-space of Step~1. By Theorem~\ref{thm:SW} and partial summation the sum over $p^*$ is $\ll_A(1+|s|)L^{-A}$ times its trivial bound. Bounding the $j$-sum trivially and counting $O(R^2)$ characters, the block is $\ll_ATR^{3/2}L^{C-A}/(g\delta ae)$. Choose $C$, then $B$, then $A$. The non-principal part is $o(TL^2)$, and the theorem follows. The $o(1)$ is ineffective.
\qed

\subsection{Remarks}
\emph{The cut.} The cut must be on the size of $n$, at a power strictly between $T^\theta$ and $T^{1/2}$. A cut at $Y=T^\theta$ does not work: take $h=p$ with $\frac34Y<p<Y$ and $n=2p$. Then $n>Y$, but after removing the shared prime $p$ the numerator is $2$, and no estimate for character sums over primes applies. With $\beta=\theta$ the third term of~\eqref{eq:LSblock} does not save when $P_0\asymp R\asymp Y$.

\emph{GRH.} Under GRH the non-principal part is $O(T^{1/2+\varepsilon}Y)$, by bounding the character sums over primes pointwise on $\re s=\frac12+\varepsilon$, and Steps~1, 2, 8 and~9 are not needed. This was the first proof of the case $d=1$. The large sieve and Siegel--Walfisz are the usual substitute for GRH; Benli, Elma and Ng use Siegel--Walfisz in the same way for their error term.

\emph{The mollifier.} The argument does not extend to $\theta>\frac12$: the fourth term of~\eqref{eq:LSblock} is then larger than $T$, and only square-root cancellation over the characters themselves, beyond what GRH gives for each character, would help. We return to this in Section~\ref{sec:remarks}.

\section{Higher moments and the proof of Theorem~\ref{thm:main}}\label{sec:higher}

\subsection{The idea}
The certificate of Section~\ref{sec:second} uses two facts about a multiple zero: the detectors vanish there, and $\re(J+a)\le0$ there. It uses nothing about the simple zeros. Short Dirichlet polynomials $U$, of length $T^\theta$, carry no pointwise information at multiple zeros. They still help, for the following reason. Since the detectors vanish off the simple zeros, they only see $\1_{\mathrm s}U$, and projecting $\1_{\mathrm s}U$ onto the detectors gives a lower bound for $\E\1_{\mathrm s}|U|^2$. On the other hand,
\[
  \E\1_{\mathrm s}|U|^2\ \le\ \E\1_{\mathrm s}+\E\bigl(|U|^2-1\bigr)_+ ,
\]
and the positive part is small if $|U|^2$ is concentrated near $1$. The positive part is controlled by the second and fourth moments of $U$, which are known for short $U$. Comparing the two bounds gives a lower bound for $p=\E\1_{\mathrm s}$. This is the clipping certificate.

The witness cannot be put inside the positive part: its fourth moments are not known. The way around this is to arrange the quadratic form so that the witness occurs only linearly, and to take the positive part only of a form in the short polynomials. This is Lemma~\ref{lem:combined}. Finally, all the positive semidefinite constraints that the available moments impose can be collected into one semidefinite program (Proposition~\ref{prop:relax}), whose dual solutions are the certificates.

The moments needed are fourth moments of short polynomials, mixed moments of the detectors against short polynomials, and, for Theorem~\ref{thm:main}, mixed moments of the detectors against products of two short polynomials. The last have length up to $T^{2\theta}$, beyond the range of the theorem of Benli, Elma and Ng; they are given by Theorem~\ref{thm:longmixed}.

\subsection{Letters}
For $d\ge1$ and a symmetric polynomial $k$ in $d$ variables, the \emph{letter} $J_d[k]$ is
\[
  J_d[k](s)=\sum_{\substack{p_1<\dots<p_d\\u_{p_1}+\dots+u_{p_d}\le\theta}}u_{p_1}\cdots u_{p_d}\,k(u_{p_1},\dots,u_{p_d})\,(p_1\cdots p_d)^{-s}.
\]
The constant $1$ is the letter of degree $0$. We use the sixteen letters of Table~\ref{tab:letters}, where
\[
  h(u)=1+\frac u{3\theta}\Bigl(1-\frac u\theta\Bigr).
\]

\begin{table}[ht]
\centering
\footnotesize
\setlength{\tabcolsep}{3pt}
\begin{tabular}{cccccccc}
\toprule
$\mathsf P$&$\mathsf S$&$\mathsf T$&$\mathsf V$&$\mathsf X$&$\mathsf Y$&$\mathsf W$&$\mathsf Z$\\
$J_1[h]$&$J_2[1]$&$J_3[1]$&$J_2[\frac{u_1+u_2}2]$&$J_1[u]$&$J_1[1]$&$J_2[u_1u_2]$&$J_2[\frac{u_1^2+u_2^2}2]$\\
\midrule
$\mathsf A$&$\mathsf F$&$\mathsf U$&$\mathsf K$&$\mathsf L$&$\mathsf D$&$\mathsf M$&$\mathsf N$\\
$J_3[\frac{u_1+u_2+u_3}3]$&$J_3[\frac{u_1^2+u_2^2+u_3^2}3]$&$J_3[\frac{u_1u_2+u_1u_3+u_2u_3}3]$&$J_1[u^3]$&$J_1[u^4]$&$J_4[1]$&$J_2[\frac{u_1u_2(u_1+u_2)}2]$&$J_2[\frac{u_1^3+u_2^3}2]$\\
\bottomrule
\end{tabular}
\caption{The letters.}\label{tab:letters}
\end{table}

A letter of degree $d$ has the kernel $k\cdot\1[u_1+\dots+u_d\le\theta]$ in the sense of Theorem~\ref{thm:longmixed}. More generally, we say that a Dirichlet polynomial with real coefficients $a_n$, supported on $n\le T^{1-\varepsilon}$, has the \emph{kernel} $\mathcal K=(\mathcal K^{(r)})_{r\le R}$ if $|a_n|\ll u(n)$, $a_n=0$ when $\Omega(n)>R$, and $a_n=u(n)\mathcal K^{(r)}(u_{p_1},\dots,u_{p_r})$ for squarefree $n=p_1\cdots p_r$; here $u(n)=\prod_{p^j\parallel n}u_p^j$, and each $\mathcal K^{(r)}$ is bounded, symmetric and continuous off finitely many hyperplanes. If $\mathcal A_1,\mathcal A_2$ have kernels $\mathcal K_1,\mathcal K_2$, then $\mathcal A_1\mathcal A_2$ has the kernel
\begin{equation}\label{eq:prodkernel}
  (\mathcal K_1\star\mathcal K_2)^{(r)}(\mathbf u)=\sum_{S\subseteq\{1,\dots,r\}}\mathcal K_1^{(|S|)}(\mathbf u|_S)\,\mathcal K_2^{(r-|S|)}(\mathbf u|_{S^c}).
\end{equation}

\begin{proposition}\label{prop:pairing}
Assume RH. If $\mathcal A$ and $\mathcal C$ have kernels $\mathcal K_A$ and $\mathcal K_C$, then $\EW\mathcal A\overline{\mathcal C}\to\langle\mathcal K_A,\mathcal K_C\rangle$, where, with $d\mu_r(\mathbf u)=\prod_{i\le r}u_i\,du_i$,
\begin{multline*}
  \langle\mathcal K_A,\mathcal K_C\rangle=\sum_{r\ge0}\frac1{r!}\int\mathcal K_A^{(r)}\mathcal K_C^{(r)}\,d\mu_r\\-\sum_{r\ge0}\frac1{r!}\iint\Bigl[\mathcal K_A^{(r+1)}(\mathbf u,w)\mathcal K_C^{(r)}(\mathbf u)+\mathcal K_A^{(r)}(\mathbf u)\mathcal K_C^{(r+1)}(\mathbf u,w)\Bigr]w\,dw\,d\mu_r(\mathbf u).
\end{multline*}
\end{proposition}

\begin{proof}
Apply Lemma~\ref{lem:errorsum}. The error is $o(1)$ after normalisation, since both lengths are at most $T^{1-\varepsilon}$. In the diagonal, the non-squarefree $n$ contribute $O(L^{-2})$, because $\sum_pu_p^4/p^2\ll L^{-4}$ and $\sum_{\Omega(n)\le R}u(n)^2/n\ll1$. For squarefree $n=p_1\cdots p_r$ we have $a_nc_n/n=\prod_i(\log p_i/Lp_i)\,u_{p_i}\cdot\mathcal K_A^{(r)}\mathcal K_C^{(r)}$, and~\eqref{eq:mertens} gives the first sum. In the off-diagonal terms, $\ell=p^j$ with $j\ge2$, $\ell\mid n$, and non-squarefree $n$ all contribute $O(L^{-1})$. In the remaining terms $\ell=q$ is a prime not dividing the squarefree $n$, and $\Lambda(q)/(Lq)=(\log q/Lq)$ contributes the variable $w=u_q$ with the weight $u_q$ coming from $a_{qn}$ or $c_{qn}$. Equation~\eqref{eq:mertens} in dimension $r+1$ gives the second sum.
\end{proof}

All entries of the second and fourth moments of letters are therefore integrals of polynomials over polytopes defined by the constraints $\sum_{i\in E}u_i\le\theta$, with rational coefficients. We evaluated them exactly. The covariances of the witness $q=J+a$ with letters and with products of letters are given by the same proposition, since $J$ has length $T^\psi<T$; by Lemma~\ref{lem:Jmom} only its prime part matters.

\emph{Mixed moments.} The mixed moment of a detector with a letter of degree $d$ is $\frac1{d!}\int k\,D_d$ over the simplex $\{u_1+\dots+u_d\le\theta\}$, by Theorem~\ref{thm:longmixed} (in fact by its Step~2 alone). A product of two letters of degrees $d_1,d_2$ is, apart from its terms with a repeated prime, a sum of polynomials as in Theorem~\ref{thm:longmixed}, of degree $d_1+d_2$, with $\psi=2\theta<1$ and the kernel~\eqref{eq:prodkernel}. The terms with a repeated prime have $\EW$-norm $O(L^{-1/2})$, so their mixed moments with detectors tend to $0$ by Cauchy--Schwarz. So Theorem~\ref{thm:longmixed} gives all mixed moments of detectors against products of two letters.

\emph{All tapers at once.} For fixed $\theta$ and $g$, the detector Gram form of Lemma~\ref{lem:BHBW} is, in terms of the tapers, $(f_1,f_2)\mapsto A_{f_1}A_{f_2}+V_{f_1f_2}$, where $A_f=\frac{f(0)}2+\int_0^\theta f$ and $V_{f_1f_2}=\frac1{12}\bigl(f_1(0)f_2(0)+\int_0^\theta f_1'f_2'\bigr)$. Every mixed moment computed above is a functional of $f$ of the form $\alpha f(0)+\int_0^\theta f\varrho$, with $\varrho$ piecewise polynomial. Its representer for the Gram form solves a second order differential equation with piecewise polynomial right side, and is computed exactly. So the norm of the projection of any finite set of these functionals onto the closure of the taper space is an exact rational number. Every certificate below is continuous in the projected data, so it suffices to compute it for the closure; finite sets of polynomial tapers approximate the closure as closely as we like. We call this \emph{using all detectors}.

\subsection{The combined certificate}
In this subsection $\E$ is any probability measure on a finite multiset of zeros with ordinates in $[T,2T]$, $p=\E\1_{\mathrm s}$, and $q=J+a$. By Lemma~\ref{lem:local},
\begin{equation}\label{eq:qsign}
  \re q\le\frac{c_0}L\ \text{at multiple zeros},\qquad \re q\le\frac12+\frac{c_0}L\ \text{at simple zeros}.
\end{equation}
For a column vector $v$ of functions on the zeros put $v^*=\bar v^{\mathsf T}$. A real symmetric matrix is PSD if it is positive semidefinite.

\begin{lemma}\label{lem:combined}
Let $H_1,\dots,H_d$ be functions on the zeros, $v_{\mathrm s}=(1,H_1,\dots,H_d)^{\mathsf T}$ and $v=(v_{\mathrm s},q)$. Let $F=(F_1,\dots,F_r)^{\mathsf T}$ be functions vanishing at all multiple zeros, with $\mathcal D:=\re\E[FF^*]$ positive definite. Put $O=\re\E[vv^*]$, $m=\re\E[Fv^*]$ and $M=m^{\mathsf T}\mathcal D^{-1}m$. Let $A,B$ be real PSD matrices, $R_{\mathrm s}$ a real symmetric matrix and $\lambda\ge0$ such that, identically in $v$ with $v_0=1$,
\begin{equation}\label{eq:combid}
  v^*(A-B)v=1+R-2\lambda\re q,\qquad R:=v_{\mathrm s}^*R_{\mathrm s}v_{\mathrm s}.
\end{equation}
Then
\[
  p\ \ge\ \tr(MA)-\tr(OB)+2\lambda\,\E\re q-\tfrac12\E R-\tfrac12\sqrt{\E R^2}-\frac{2\lambda c_0}{L}.
\]
\end{lemma}

\begin{proof}
Write $A=\sum_ka_ka_k^{\mathsf T}$ with real vectors $a_k$, and put $X_k=a_k^{\mathsf T}v$. For real $\beta$, since $F=\1_{\mathrm s}F$,
\[
  0\le\E\bigl|\1_{\mathrm s}X_k-\beta^{\mathsf T}F\bigr|^2=\E\1_{\mathrm s}|X_k|^2-2\beta^{\mathsf T}ma_k+\beta^{\mathsf T}\mathcal D\beta .
\]
With $\beta=\mathcal D^{-1}ma_k$ this gives $\E\1_{\mathrm s}|X_k|^2\ge a_k^{\mathsf T}Ma_k$, and so $\E\1_{\mathrm s}v^*Av\ge\tr(MA)$. Since $v^*Bv\ge0$ pointwise, $\E\1_{\mathrm s}v^*Bv\le\tr(OB)$. By~\eqref{eq:combid} and~\eqref{eq:qsign},
\[
  \1_{\mathrm s}v^*(A-B)v=\1_{\mathrm s}+\1_{\mathrm s}R-2\lambda\re q+2\lambda(1-\1_{\mathrm s})\re q\le\1_{\mathrm s}+R_+-2\lambda\re q+\frac{2\lambda c_0}L .
\]
Finally $R_+\le(R+r)^2/4r$ for every $r>0$, and $r=\sqrt{\E R^2}$ gives $\E R_+\le\frac12(\E R+\sqrt{\E R^2})$. Take expectations and combine.
\end{proof}

The witness enters~\eqref{eq:combid} only linearly, through the cross term with the constant. So $\E R^2$ is a fourth moment of short polynomials, and the unknown fourth moments of $J$ are never needed. With $\lambda=0$, $q$ absent, $A=UU^{\mathsf T}$, $B=0$ and $R=|U|^2-1$, the lemma is the clipping certificate described above. With $R=0$, it reduces to the second-order certificate.

\subsection{The semidefinite relaxation}
Let $v_0=1$, $v_i=H_i$ ($1\le i\le d$), and
\[
  w=\bigl(v_iv_j\bigr)_{0\le j\le i\le d}\oplus(q),\qquad Y=\re\E[\1_{\mathrm s}ww^*],\qquad O_w=\re\E[ww^*].
\]
Then $w$ contains $v_{\mathrm s}=(v_iv_0)_i$ and $q$, and $Y_{00}=p$. Let $M_w=m_w^{\mathsf T}\mathcal D^{-1}m_w$ with $m_w=\re\E[Fw^*]$, and for a set $\mathrm{sel}$ of coordinates of $w$ let $Y_{\mathrm{sel}}$ and $M_{\mathrm{sel}}$ be the corresponding principal submatrices.

\begin{proposition}\label{prop:relax}
The matrix $Y$ satisfies:
\begin{enumerate}[label=\textup{(\roman*)},leftmargin=*]
\item $0\preceq Y\preceq O_w$ and $Y_{00}=p$;
\item $Y_{\mathrm{sel}}\succeq M_{\mathrm{sel}}$ for every $\mathrm{sel}$;
\item $\E\re q-c_0/L\le Y_{0q}\le p/2+c_0/L$;
\item with $X_{ij}=v_i\bar v_j$, the Gram matrices of the real functions $\re X_{ij}$ $(i\ge j)$ and $\im X_{ij}$ $(i>j)$ with respect to the measures $\1_{\mathrm s}\,d\E$ and $(1-\1_{\mathrm s})\,d\E$ are PSD, and their entries are fixed linear combinations of the entries of $Y$, respectively of $O_w-Y$.
\end{enumerate}
\end{proposition}

\begin{proof}
(i) $Y$ and $O_w-Y$ are real parts of Hermitian PSD matrices. (ii) is the first step of the proof of Lemma~\ref{lem:combined}. (iii) follows from~\eqref{eq:qsign}, since $Y_{0q}=\E\1_{\mathrm s}\re q$. For (iv), $\re X_{ij}\re X_{kl}=\frac12\re\bigl((v_iv_l)\overline{(v_jv_k)}+(v_iv_k)\overline{(v_jv_l)}\bigr)$, and similarly $\im X_{ij}\im X_{kl}=\frac12\re\bigl((v_iv_l)\overline{(v_jv_k)}-(v_iv_k)\overline{(v_jv_l)}\bigr)$.
\end{proof}

Thus $p$ is at least the minimum of $Y_{00}$ over all real symmetric $Y$ satisfying (i)--(iv). This is a semidefinite program, and its data $O_w$, $M_{\mathrm{sel}}$ and $\E\re q$ converge as $T\to\infty$. A lower bound is certified by a dual solution: PSD multipliers $Z_{\mathrm{sel}}$ for the constraints (ii), $Z'$ for $Y\preceq O_w$, PSD multipliers for (iv), and $\lambda,b\ge0$ for (iii). Pairing each constraint with its multiplier gives an identity
\[
  Y_{00}=c+\sum_{\text{constraints}}\langle\text{multiplier},\text{constraint}\rangle+\sum_{i,j}\kappa_{ij}Y_{ij},
\]
valid for all symmetric $Y$, where $c$ is an affine function of the data and each pairing is nonnegative on the feasible set. For an exact dual solution $\kappa=0$. Our multipliers are rational, and $\kappa$ is small but not zero. Since $0\preceq Y\preceq O_w$, we have $|Y_{ij}|\le\sqrt{(O_w)_{ii}(O_w)_{jj}}$, and hence
\begin{equation}\label{eq:dualbound}
  p\ \ge\ c-\varepsilon-O(L^{-1}),\qquad \varepsilon:=\sum_{i,j}|\kappa_{ij}|\sqrt{(O_w)_{ii}(O_w)_{jj}} .
\end{equation}
The dual of the program with (iv) omitted and with a rank-one multiplier for (ii) is Lemma~\ref{lem:combined}.

\subsection{Short polynomials only}
The first theorem uses no mixed moments with long composite numerators.

\begin{theorem}\label{thm:fourth}
Assume RH. Then $\liminf_{T\to\infty}N_s(T)/N(T)\ge0.7094399755$.
\end{theorem}

\begin{proof}
Work at fixed $\theta<\frac12<\psi<1$, let $T\to\infty$, and then let $(\theta,\psi)\to(\frac12,1)$, as in~\eqref{eq:supcert}. Take the kernel
\begin{multline*}
  g(z)=1-\tfrac{31341}{200000}\Leg_2(2z)+\tfrac{26639}{500000}\Leg_4(2z)-\tfrac{3747}{200000}\Leg_6(2z)\\
  +\tfrac{363}{31250}\Leg_8(2z)-\tfrac1{125}\Leg_{10}(2z)+\tfrac3{500}\Leg_{12}(2z),
\end{multline*}
the seven letters $\mathsf{P,S,T,V,X,Y,A}$, and all detectors. Impose (ii) of Proposition~\ref{prop:relax} only for $\mathrm{sel}=(1,\mathsf P,\mathsf S,\mathsf T,\mathsf V,\mathsf X,\mathsf Y,\mathsf A,q)$. Then the data of the program are given by Lemma~\ref{lem:BHBW}, Lemma~\ref{lem:Jmom}, Proposition~\ref{prop:pairing}, Corollary~\ref{cor:C} and the mixed moments of detectors with single letters. The products of letters enter only through $O_w$, (i) and (iv). For the rational dual solution described in Appendix~\ref{app:data}, the limiting values at $(\theta,\psi)=(\frac12,1)$ are
\[
  c=0.70943998950468955\ldots,\qquad \varepsilon=1.40020122558514\ldots\times10^{-8},
\]
so that $c-\varepsilon=0.70943997550267729\ldots$. By~\eqref{eq:dualbound} and Lemma~\ref{lem:unweight} this proves the theorem.
\end{proof}

Without the semidefinite constraints (iv), Lemma~\ref{lem:combined} with nine letters and fourteen detectors already gives $0.7093885952$. Most of the gain over Theorem~\ref{thm:second} comes from cancelling the witness before the positive part is taken.

\subsection{Proof of Theorem~\ref{thm:main}}
Now we also use the mixed moments of the detectors against products of letters, that is, (ii) for all coordinates of $w$. With the seven letters of Theorem~\ref{thm:fourth} and a similar kernel this gives $0.7126403094$. To do better we use all sixteen letters of Table~\ref{tab:letters}, but compress them, since the program with sixteen letters and all their products is large. A rational $8\times17$ matrix $R$, with first row $(1,0,\dots,0)$, maps the constant and the sixteen letters to the constant and seven \emph{probes} $v_1,\dots,v_7$, which are fixed linear combinations of letters. Their products are linear combinations of products of letters, through the symmetric square of $R$. With $\Lambda=\mathrm{Sym}^2R\oplus1$, the data of the program for the probes are exactly $\Lambda O_w\Lambda^{\mathsf T}$ and $\Lambda M_w\Lambda^{\mathsf T}$, where $O_w$ and $M_w$ are the data for the sixteen letters.

All columns of $M_w$ are used. The longest products, $\mathsf D\cdot\mathsf D$, have eight prime factors and length $T^{2\theta}$; their mixed moments with detectors come from Theorem~\ref{thm:longmixed} with $d=8$ and $\psi=2\theta$. The kernel is
\begin{multline*}
  g(z)=1-\tfrac{4639}{25000}\Leg_2(2z)+\tfrac{6847}{100000}\Leg_4(2z)-\tfrac{2811}{100000}\Leg_6(2z)\\
  +\tfrac{1861}{100000}\Leg_8(2z)-\tfrac{67}{5000}\Leg_{10}(2z)+\tfrac{13}{1250}\Leg_{12}(2z).
\end{multline*}
The dual multipliers are those of a rational dual solution for the seven letters $\mathsf{P,S,T,V,X,Y,A}$ themselves; only $R$ was optimised. At $(\theta,\psi)=(\frac12,1)$ the certificate has
\[
  c=0.71274069505648300\ldots,\qquad \varepsilon=2.66376033908722\ldots\times10^{-8},
\]
so that $c-\varepsilon=0.71274066841887961\ldots$. By~\eqref{eq:dualbound}, letting $(\theta,\psi)\to(\frac12,1)$ after $T\to\infty$, and Lemma~\ref{lem:unweight}, this proves Theorem~\ref{thm:main}. \qed

The products of letters are what moves the bound from $0.7094$ to $0.7127$. They enter because a product of two letters is again a short polynomial, of length $T^{2\theta}<T$, whose pairing with the detectors we can compute. In the language of the introduction, they allow the certificate to use quadratic functions of the short polynomials as additional detectors of simplicity, constrained by the fourth moments.

\section{Distinct zeros}\label{sec:distinct}

A lower bound for $N_s/N$ does not by itself give a lower bound for $N_d/N$ better than $N_s/N$, since the multiple zeros might have large multiplicity. The certificate of Section~\ref{sec:second} gives more, because its residual is large at zeros of high multiplicity. In addition, double zeros cannot be packed densely: three nearby double zeros make the kernel sum in~\eqref{eq:local} larger than the multiplicity. We use both facts.

\subsection{The residual at a multiple zero}
Let $K$ be the function attached to the kernel $g_*$ of~\eqref{eq:gstar} as in Section~\ref{sec:kernels}, and fix the optimal detector coefficients and the witness coefficient $t=\tau_*$ of the certificate~\eqref{eq:p0} at $(\theta,\psi)=(\frac12,1)$. For $\theta<\frac12<\psi<1$ consider the residual
\[
  \rho_0=1-F-t(J+a),
\]
where $F$ is the corresponding detector and $J$, $a$ are built from $K_\psi$. Then $\EW|\rho_0|^2\to1-\kappa(\theta,\psi)$, with $\kappa(\theta,\psi)\to p_0$ as $(\theta,\psi)\to(\frac12,1)$. Write
\[
  D(\rho)=\sum_{\gamma'}K_\psi\Bigl(\frac{(\gamma'-\gamma)L}{2\pi}\Bigr)\ \ge\ m_\rho .
\]
At a multiple zero, $F=0$ and, by~\eqref{eq:local}, $\re(J+a)\le1-D/2+c_0/L$. Hence $\re\rho_0\ge1+t(D-2)/2-c_0t/L$, and since $(1+x)^2\ge1+2x$,
\begin{equation}\label{eq:resid}
  |\rho_0|^2\ \ge\ 1+t(m_\rho-2)+t(D-m_\rho)-O(L^{-1})\qquad(m_\rho\ge2).
\end{equation}
The term $t(m_\rho-2)$ rewards high multiplicity, and $t(D-m_\rho)$ rewards nearby zeros.

\subsection{Packing of double zeros}
\begin{lemma}\label{lem:gap}
Let $b_1=1048576/21217072921875$ and $A=16$. For all $x,y\ge0$,
\[
  \max\bigl(K(x),K(y),K(x+y)\bigr)\ \ge\ b_1\Bigl(1-\frac{x+y}A\Bigr).
\]
\end{lemma}

\begin{proof}
For $x+y\ge A$ there is nothing to prove. On the triangle $x,y\ge0$, $x+y\le A$, the function $K(x)=\check g_*(x)^2$ is an explicit combination of $\sin\pi x$, $\cos\pi x$ and negative powers of $x$. We subdivided the triangle into cells and bounded $K$ from below on each cell in exact rational interval arithmetic.
\end{proof}

The only reason for the lemma is that $K$ has no common zero at $x$, $y$ and $x+y$: two consecutive gaps and their sum cannot all be zeros of $K$. The constant $b_1\approx4.9\times10^{-8}$ is small because $K$ is small beyond the first few zeros of $\check g_*$.

Let $p_m=\EW\1_{m\ge2}$, $h_3=\EW\1_{m\ge3}$ and $d_2=(p_m-h_3)/2$. So $d_2N_W$ is the weighted number of distinct double zeros, which we call \emph{double sites}.

\begin{lemma}\label{lem:packing}
Assume RH. For $\theta<\frac12<\psi<1$,
\[
  \EW\bigl[\1_{m\ge2}(D-m)\bigr]\ \ge\ 4b_1\Bigl(d_2-\frac{2\psi}A\Bigr)-o(1).
\]
\end{lemma}

\begin{proof}
Partition $[T,2T]$ into intervals of length $\eta T$, on each of which $W(\gamma/T)$ is constant up to $O(\eta)$, and order the double sites in each interval. Measure distances in units of $2\pi/L$, the mean spacing. Let $x,y$ be the gaps between three consecutive double sites. By Lemma~\ref{lem:gap} applied to $(\psi x,\psi y)$,
\[
  \tfrac12\bigl(K_\psi(x)+K_\psi(y)\bigr)+K_\psi(x+y)\ \ge\ \tfrac12\max\bigl(K_\psi(x),K_\psi(y),K_\psi(x+y)\bigr)\ \ge\ \tfrac12b_1\Bigl(1-\frac{\psi(x+y)}A\Bigr).
\]
Sum over all consecutive triples. On the left, each pair of adjacent double sites occurs in at most two triples with weight $\frac12$, and each pair at distance two in at most one triple with weight $1$. So the left side is at most $\sum K_\psi(\text{distance})$ over unordered pairs of double sites. Each such pair contributes $8$ times its $K_\psi$-value to $\sum_\rho\1_{m\ge2}(D-m)$, a sum over zeros counted with multiplicity: each site has multiplicity two and sees the other with multiplicity two. On the right, the two-gaps $x+y$ cover each point at most twice, and the total length is the number of zeros up to $o(N_W)$. The boundary losses are $O(1)$ per interval. Let $T\to\infty$ and then $\eta\to0$.
\end{proof}

No hypothesis on the spacing of the zeros is used.

\subsection{Proof of Theorem~\ref{thm:distinct}}
Let $\mathrm d=\EW[1/m_\rho]$. For $m\ge2$ and $t\ge\frac16$,
\[
  1+t(m-2)\ \ge\ 2\Bigl(1-\frac1m\Bigr)+\Bigl(t-\frac13\Bigr)\1_{m\ge3},
\]
with equality at $m=2$ and $m=3$. Also $1-\mathrm d=\EW[1-1/m]\le d_2+h_3$. Take expectations in~\eqref{eq:resid} and use Lemma~\ref{lem:packing}. With $C_1=2tb_1$,
\[
  1-\kappa(\theta,\psi)\ \ge\ 2(1-\mathrm d)+\Bigl(t-\frac13\Bigr)h_3+2C_1\Bigl(1-\mathrm d-h_3-\frac{2\psi}A\Bigr)-o(1).
\]
Since $t=0.7387\ldots$, the coefficient $t-\frac13-2C_1$ of $h_3$ is positive, and
\[
  \liminf_{T\to\infty}\mathrm d\ \ge\ \frac12\Bigl(1+\frac{\kappa(\theta,\psi)+C_1(1-4\psi/A)}{1+C_1}\Bigr).
\]
Let $(\theta,\psi)\to(\frac12,1)$ and apply Lemma~\ref{lem:unweight} with $\phi=1/m_\rho$. Since
\[
  S:=\frac{p_0+C_1(1-4/A)}{1+C_1}=0.70736460920095\ldots,
\]
we obtain $\liminf N_d(T)/N(T)\ge(1+S)/2=0.85368230460047\ldots$. \qed

The same argument with $1+t(m-2)\ge1+t\1_{m\ge3}$ gives $\liminf N_s/N\ge S$, which is weaker than Theorem~\ref{thm:main}. Without packing, the certificate of Theorem~\ref{thm:second} gives $N_d/N\ge(1+\kappa)/2=0.8536822884\ldots$, since its witness coefficient exceeds $\frac13$, and the kernel $g_*$ gives $(1+p_0)/2=0.8536823030\ldots$. The packing adds about $1.5\times10^{-9}$. The certificates of Section~\ref{sec:higher} could be used in the same way, but their residuals at multiple zeros are harder to control, and we have not done so.

\section{Remarks and open problems}\label{sec:remarks}

\subsection{The second-order ceiling}
Theorem~\ref{thm:second} is close to the best possible for its method.

\begin{proposition}\label{prop:ceiling}
At $(\theta,\psi)=(\frac12,1)$, for every kernel $k$ that is a finite sum of autocorrelations $g*g$ of real even functions supported in $[-\frac12,\frac12]$, and every finite set of tapers, the value $\Ksq/(1+\Ksq)$ of Proposition~\ref{prop:optim} is less than $0.7073646063$.
\end{proposition}

Kernels with $\psi<1$ are included, since $k_\psi$ is the autocorrelation of a function supported in $[-\frac\psi2,\frac\psi2]$.

\begin{proof}[Sketch of proof]
Fix $k=g*g$ with $\int g=1$. By Proposition~\ref{prop:optim}, the optimum over all tapers is $\mathcal R(\mathcal A,\mathcal A)+\gap_+^2/D_J$, where $\mathcal R$ is the inverse of the Gram form of the tapers described in Section~\ref{sec:higher} (\emph{All tapers at once}), and $\mathcal A$, $\mathfrak c_k$ are the representers of the detector means and of the covariances of Corollary~\ref{cor:C}. Using $2\int_0^1k=(\int g)^2=1$ we write
\[
  \gap=G(k):=2\int_0^1k-\frac{k(0)}2-\int_0^1uk-\mathcal R(\mathcal A,\mathfrak c_k),\qquad D_J=D(k,k),
\]
\[
  D(k,l):=\int_0^1ukl-\int_0^1uk\int_0^1ul-\mathcal R(\mathfrak c_k,\mathfrak c_l).
\]
Now $G$ is linear and $D$ is quadratic in $k$, so these expressions extend to finite sums of autocorrelations. The form $D$ is positive semidefinite on them, being the variance of the witness after projection onto the detectors. Scaling $k$, the increment $\gap_+^2/D_J$ is the supremum of $2G(k)-D(k,k)$ over the cone of finite sums of autocorrelations. If $s$ is a trial function with $G(k)\le D(s,k)$ for all $k$ in the cone, then this supremum is at most $D(s,s)$, because $D(k-s,k-s)\ge0$. We take for $s$ a multiple of an optimised autocorrelation. For $k=g*g$, the condition $G(k)\le D(s,k)$ says that an explicit integral operator $T_E$ on $L^2[0,1]$, with kernel $E(|x-y|)$ for an explicit piecewise polynomial $E\ge0$, satisfies $T_E\le I$. We verify this by a weighted Schur test: $T_Eh\le h$ for an explicit positive polynomial $h$ of degree $22$, checked by subdivision and Bernstein coefficients in exact arithmetic. The resulting rational upper bound is $0.70736460621287\ldots$.
\end{proof}

The kernel $g_*$ gives $p_0=0.7073646060\ldots$, within $2\times10^{-10}$ of the ceiling, and the trial kernel used in the proof gives $0.7073646062056\ldots$, within $10^{-11}$. So the second-order method is exhausted, and higher moments are necessary to go further.

\emph{Sharpness for fixed moments.} For a fixed finite set of tests, the bound $m^{\mathsf T}G^{-1}m$ is the best that the moments $(m,G)$ alone can give. Indeed, put $\kappa=m^{\mathsf T}G^{-1}m$ and consider a random vector $X$ that vanishes with probability $1-\kappa$ and, with probability $\kappa$, has mean $m/\kappa$ and second moment matrix $G/\kappa$. This is possible because $G-mm^{\mathsf T}/\kappa\succeq0$. Then $X$ has the prescribed moments and vanishes on a set of measure $1-\kappa$. Any improvement must therefore use more information, as in Section~\ref{sec:higher}.

\emph{The long mixed moment is needed.} One might try to avoid Theorem~\ref{thm:longmixed} by bounding the covariance of $F$ with the long part of $J$ by Cauchy--Schwarz. In the proof sketch of Section~\ref{sec:sketch}, the long part $J_{\mathrm{long}}$ (primes in $(T^{1/2},T)$) has $\E|J_{\mathrm{long}}|^2\to\int_{1/2}^1u(1-u)^2du=\frac5{192}$, so the covariance would be known only up to $\pm\sqrt{V\cdot5/192}=\pm0.0828\ldots$. This is five times the gap $\frac{31}{1920}=0.0161\ldots$, and the witness would give nothing.

\subsection{Longer mollifiers and pair correlation}
Two limits are built into the method. The mollifier has length $T^\theta$ with $\theta<\frac12$, and the witness has length $T^\psi$ with $\psi<1$.

Let $H(\theta)$ denote the hypothesis that Lemma~\ref{lem:BHBW} and Theorem~\ref{thm:longmixed} hold with the given $\theta>\frac12$. Under $H(\theta)$ all of our certificates improve, and the mixed moments with long numerators become more valuable, since the detectors see more primes. The obstruction is the one already met by Conrey, Ghosh and Gonek. In the proof of Theorem~\ref{thm:longmixed} the moduli are the mollifier variables, and for moduli beyond $T^{1/2}$ the term $R\sqrt T$ of~\eqref{eq:LSblock} exceeds the main term. Under GRH each character sum is bounded individually, but the number of characters is the problem, so GRH does not help either. What would be needed is square-root cancellation in the sum over characters of the products of the two character sums, beyond the large sieve.

Witnesses with $\psi>1$ would require the pair correlation function $F(\alpha)$ for $\alpha>1$. This is Montgomery's pair correlation conjecture, which implies that almost all zeros are simple~\cite{Mon73}. The pointwise use of the explicit formula does not avoid this: the witness of length $T^\psi$ has known moments only for $\psi<1$.

Nothing here is unconditional. The unconditional results for simple zeros, starting with Conrey's $\frac25$~\cite{Con89}, use mollified moments off the critical line, and we do not see how to combine them with the pointwise inequality, which needs RH.

\subsection{Dirichlet $L$-functions}\label{sec:family}
The local argument applies to any family where the pair correlation is known in a range. For primitive Dirichlet characters, Chandee, Lee, Liu and Radziwi\l\l~\cite{CLLR14} computed the pair correlation of the low-lying zeros averaged over the family, under GRH, for Fourier support in $[-2,2]$.
Let $W\ge0$ be smooth and supported in $(1,2)$, $w_q=W(q/Q)/\varphi(q)$, and let $a_\gamma=\widehat\Phi(i\gamma)\ge0$ be the weight of~\cite{CLLR14}. Put
\[
  N_\Phi(Q)=\sum_qw_q\sideset{}{^*}\sum_{\chi\bmod q}\ \sum_{\gamma_\chi}a_{\gamma_\chi}^2,
\]
where the inner sum runs over the zeros $\frac12+i\gamma_\chi$ of $L(s,\chi)$, counted with multiplicity. Let $N_{\Phi,s}(Q)$ be the same sum over simple zeros, and $N_{\Phi,d}(Q)$ the sum over distinct zeros, each counted once. Proposition~\ref{prop:family} states that, under GRH,
\begin{align*}
  \liminf_{Q\to\infty}\frac{N_{\Phi,s}(Q)}{N_\Phi(Q)}&\ge\tfrac{29246037121493404426049952732657534866284336504156446769}{31276240776124426360800941323301439283142863495843553231}\\&=0.935087989980\ldots,\\
  \liminf_{Q\to\infty}\frac{N_{\Phi,d}(Q)}{N_\Phi(Q)}&\ge\tfrac{30261138948808915393425447027979487074713600000000000000}{31276240776124426360800941323301439283142863495843553231}\\&=0.967543994990\ldots.
\end{align*}
For comparison, Chirre, Gon\c{c}alves and de Laat~\cite[Theorem~5]{CGdL20} obtain $0.9350$ for simple zeros.

\begin{proof}[Proof of Proposition~\ref{prop:family}]
Let $U=\log Q/2\pi$, and let $r$ be an even function with $r\ge0$, $r(0)=1$, whose Fourier transform $\hat r$ is even, piecewise linear, compactly supported, and $\le0$ outside $[-2,2]$. A zero of $L(s,\chi)$ of multiplicity $m$ satisfies $m\le\sum_{\gamma'}r(U(\gamma'-\gamma))$, the sum running over the zeros of the same $L$-function. Let $N^*_\Phi$ be $N_\Phi$ with each zero weighted by its multiplicity. Summing over the family and using the pair correlation theorem of~\cite{CLLR14} (with $r(sx)$, $s<1$, to stay inside its range, and $s\uparrow1$ at the end), we obtain
\[
  \frac{N^*_\Phi(Q)}{N_\Phi(Q)}\le B(r)+o(1),\qquad B(r)=\hat r(0)+\int_{-2}^2\min(1,|u|)\,\hat r(u)\,du .
\]
Here we used that the weight $a_\gamma^2$ may be replaced by $a_\gamma a_{\gamma'}$ at a cost $o(1)$, that the pair correlation measure is nonnegative, and that $\hat r\le0$ outside $[-2,2]$. Since $m\ge2$ off the simple zeros and $1/m\ge(3-m)/2$, this gives the Montgomery-type bounds $N_{\Phi,s}/N_\Phi\ge2-B$ and $N_{\Phi,d}/N_\Phi\ge(3-B)/2$.

To improve them we keep the off-diagonal terms from consecutive distinct zeros, as in Lemma~\ref{lem:packing}. For gaps $x,y\ge0$ in units of $1/U$ put $M_r(x,y)=\frac12(r(x)+r(y))+r(x+y)$. For the kernel $r$ described in Appendix~\ref{app:data} we have
\begin{equation}\label{eq:familypacking}
  M_r(x,y)\ \ge\ \beta\Bigl(1-\frac{x+y}A\Bigr)-\delta\qquad(x,y\ge0),
\end{equation}
with $A=\frac{13}5$, $\beta=\frac{95177}{75631303}$ and $\delta<6.63\times10^{-7}$, and
\[
  B=B(r)=\frac{2876606615715690641429157806}{2700000001435860000190897587}.
\]
The inequality~\eqref{eq:familypacking} is checked on a mesh in exact rational arithmetic, with the interpolation error bounded by $\|r''\|_\infty$; $\delta$ absorbs it. Summing~\eqref{eq:familypacking} over consecutive triples of distinct zeros of each $L(s,\chi)$, in short intervals where $a_\gamma^2$ is nearly constant, gives, with $z=\liminf N^*_\Phi/N_\Phi$ and $\mathrm d_\Phi=\liminf N_{\Phi,d}/N_\Phi$,
\[
  z+2(\beta-\delta)\,\mathrm d_\Phi-\frac{4\beta}A\ \le\ B .
\]
Combined with $\mathrm d_\Phi\ge(3-z)/2$ this gives $\mathrm d_\Phi\ge d_*$, and $N_{\Phi,s}/N_\Phi\ge2d_*-1$, where
\[
  d_*=\frac{(3-B)/2-2\beta/A}{1-\beta+\delta}.
\]
These are the stated fractions.
\end{proof}

Without packing the same $r$ gives $2-B=0.93459\ldots$, and a kernel optimised for the bound $2-B(r)$ alone gives $0.9350495091\ldots$. We have not found a way to make a mollified detector useful here. The gain comes from the larger support $[-2,2]$, which already makes the pair correlation bound strong, and from packing.

\section*{Acknowledgments}
The research framework and overall strategy of this work were set by the author. A self-developed mathematical agentic system, Dasheng, played a key role in the mathematical exploration, in particular in developing the constructions and proofs within this framework. All constants were evaluated in exact rational arithmetic; the verification code and certificate data are archived on Zenodo at \href{https://doi.org/10.5281/zenodo.24536326}{doi:10.5281/zenodo.24536326}.

\appendix
\section{Certificate data}\label{app:data}

This appendix lists the data that determine each constant. Every displayed constant is the exact value of a rational expression, truncated.

\subsection*{Data for Theorem~\ref{thm:second}}
The kernel is $g_8$ of~\eqref{eq:g8}, and the tapers are $b_0,\dots,b_{11}$ of~\eqref{eq:basis}, with $A$ and $V$ given by~\eqref{eq:basisAV}. The covariances of Corollary~\ref{cor:C} are
\begin{align*}
  C_{P_0}&=-\tfrac{9401490355717}{170311680000000}, & C_{P_1}&=\tfrac{1814282429903}{442810368000000}, & C_{P_2}&=\tfrac{5968179133261}{11070259200000000},\\
  C_{P_3}&=-\tfrac{217413232843}{2214051840000000}, & C_{P_4}&=\tfrac{1488400877}{4428103680000000}, & C_{P_5}&=\tfrac{338458349}{885620736000000},\\
  C_{P_6}&=-\tfrac{67629533}{4428103680000000}, & C_{P_7}&=-\tfrac{312934007}{420669849600000000}, & C_{P_8}&=\tfrac{147223}{12019138560000000},\\
  C_{P_9}&=-\tfrac{1181}{12019138560000000}, & C_{P_{10}}&=-\tfrac{17}{801275904000000}, & C_{P_{11}}&=-\tfrac{17}{55288037376000000}.
\end{align*}
With $\mu$, $a$ and $V_J$ as in the proof, the optimal coefficients are $G^{-1}m$ (Proposition~\ref{prop:optim}). The witness coefficient is the value of $\tau$ given in the proof, and the taper is $f=\sum_jd_jb_j$ with
\[
  d_0=-\tfrac{7430741787264374836523251391328}{6735054433292025624721100070697}=-1.1032\ldots,
\]
\[
  f(0)=\tfrac{3715370893632187418261625695664}{6735054433292025624721100070697}=0.5516\ldots .
\]

\subsection*{Data for~\eqref{eq:p0} and Theorem~\ref{thm:distinct}}
The kernel is $g_*$ of~\eqref{eq:gstar}, with all polynomial tapers vanishing at $\frac12$. The witness coefficient is
\[
  \tau_*=\tfrac{3565611778673678790271393717091077324800000000000000}{4826245010227457178772530978845119244788120650314257}.
\]
The constants in Lemmas~\ref{lem:gap} and~\ref{lem:packing} are $A=16$ and $b_1=1048576/21217072921875$. Then $C_1=2\tau_*b_1$ is $7.30245904604078\ldots\times10^{-8}$, so that $S=0.70736460920095\ldots$ and $(1+S)/2=0.85368230460047\ldots$.

\subsection*{Data for Theorems~\ref{thm:fourth} and~\ref{thm:main}}
The kernels are given in the proofs, and the letters in Table~\ref{tab:letters}. The rational dual solutions consist of a multiplier for each constraint of Proposition~\ref{prop:relax}: a rational positive semidefinite matrix for each of (i) and (iv), and a vector for each imposed instance of (ii) and (iii). For Theorem~\ref{thm:main} there is also the rational $8\times17$ matrix $R$. These data are too long to print, and are provided in the Zenodo archive (\href{https://doi.org/10.5281/zenodo.24536326}{doi:10.5281/zenodo.24536326}, files \texttt{higher\_moments/certificates/short7.json} and \texttt{refit16.json}). Given them, $c$ and $\varepsilon$ in~\eqref{eq:dualbound} are finite sums of rational numbers.

\subsection*{Data for Proposition~\ref{prop:ceiling}}
The trial function $s$ is $(1+10^{-7})G(k_*)k_*/D(k_*,k_*)$ for a rational autocorrelation $k_*$ with four Legendre modes. The Schur weight $h$ is a polynomial of degree $22$ with rational coefficients, positive on $[0,\frac12]$ and extended evenly about $\frac12$; the exact coefficients and verification script are in \texttt{higher\_moments/ceiling.py} in the Zenodo archive.

\subsection*{Data for Proposition~\ref{prop:family}}
Let $n=90$ and let $c_0,\dots,c_{500}$ be the integers in \texttt{distinct\_and\_family/local\_factor.json} in the Zenodo archive. Let $a_k=\sum_jc_jc_{j+k}$. Put $a'_k=\min(a_k,0)$ for $k\ge2n$, $a'_k=a_k$ for $1\le k<2n$, and $a'_0=a_0+2\sum_{k\ge2n}\max(a_k,0)$. The kernel is
\[
  r(x)=\Bigl(\frac{\sin(\pi x/n)}{\pi x/n}\Bigr)^2\,\frac{a'_0+2\sum_{k\ge1}a'_k\cos(2\pi kx/n)}{\bigl(\sum_jc_j\bigr)^2}.
\]
The trigonometric polynomial is $|\sum_jc_je(jx/n)|^2$ plus $2\sum_{k\ge2n}\max(a_k,0)(1-\cos(2\pi kx/n))\ge0$, so $r\ge0$, and $r(0)=1$. The Fourier transform $\hat r$ is piecewise linear with nodes at $k/n$, and it is $\le0$ outside $[-2,2]$. In~\eqref{eq:familypacking}, the inequality is checked at the points of the mesh $\frac1{16384}\mathbb Z^2$ with $x+y\le A$, and
\[
  \delta=\tfrac{21086355742804723181916840868699243368029}{31850496016938093775051919567093760000000000000}=6.6204\ldots\times10^{-7}
\]
includes the interpolation error.

\end{document}